\documentclass[reqno,12pt]{amsart}
\usepackage{graphicx} 
\usepackage{booktabs}
\usepackage[numbers,square]{natbib}
\allowdisplaybreaks

\makeatletter
\def\l@section{\@tocline{1}{0pt}{0pc}{}{}}
\def\l@subsection{\@tocline{2}{0pt}{2.3em}{}{}}
\makeatother

\usepackage{amsmath,amsfonts,amssymb,amsthm,amscd}

\usepackage{epsfig}
\usepackage{psfrag}
\usepackage{perpage}
\usepackage{url}
\usepackage{color}
\usepackage[english]{babel}
\usepackage{bbm}
\usepackage{epstopdf}
\usepackage[utf8]{inputenc}
\usepackage[T1]{fontenc}
\usepackage{microtype}
\usepackage{hyperref}
\hypersetup{colorlinks=true, linkcolor=blue, citecolor=blue,  filecolor=blue, urlcolor=blue}
\usepackage{mathabx}

\counterwithin{figure}{section}
\numberwithin{figure}{section}
\numberwithin{equation}{section}

\usepackage{fullpage}

\newcommand{\dd}{{\rm d}}
\newcommand{\I}{\mathbb{I}}

\newcommand{{\paa}[1]}{p_{00,#1}}
\newcommand{{\pab}[1]}{p_{01,#1}}
\newcommand{{\pba}[1]}{p_{10,#1}}
\newcommand{{\pbb}[1]}{p_{11,#1}}

\newcommand{\nelem}[1]{{Lemma \ref{#1}}}
\newcommand{\neprop}[1]{{Proposition \ref{#1}}}
\newcommand{\netheo}[1]{{Theorem \ref{#1}}}

\newcommand{\eee}{\mathrm{e}}

\newcommand{\abs}[1]{\left\lvert #1 \right\rvert}
\newcommand{\E}[1]{\mathbb{E}\left\{ #1\right\}}
\newcommand{\Cov}[1]{\operatorname{Cov}\left( #1\right)}
\newcommand{\covv}[1]{\operatorname{Cov}\bigl( #1\bigr)}

\newcommand{\prooftheo}[1]{\textsc{\bf Proof of Theorem} \ref{#1}:}
\newcommand{\N}{\mathbb{N}}
\newcommand{\R}{\mathbb{R}}
\newcommand{\vk}[1]{\boldsymbol{#1}}
\newcommand{\scale}[1]{\langle#1\rangle}
\newcommand{\QED}{\hfill $\Box$}

\newcommand{\pk}[1]{\mathbb{P} \left\{ #1 \right \} }

\newtheorem{theo}{Theorem}[section]
\newtheorem{lem}[theo]{Lemma}
\newtheorem{prop}[theo]{Proposition}
\newtheorem{remark}[theo]{Remark}

\title{Functional central limit theorem \\ for the subgraph count of the voter model \\ on dynamic random graphs}
\author[S. Baldassarri]{Simone Baldassarri}
\address{Gran Sasso Science Institute, Viale Francesco Crispi 7, 67100 L’Aquila, Italy}
\email{simone.baldassarri@gssi.it}
\thanks{}
\author[N. Kriukov]{Nikolai Kriukov}
\address{Korteweg-de Vries Institute for Mathematics, University of Amsterdam, Science Park 904, 1098 XH Amsterdam, The Netherlands}
\email{n.kriukov@uva.nl}

\begin{document}

\begin{abstract}
In this paper we consider two-opinion voter models on dynamic random graphs, in which the joint dynamics of opinions and graphs acts as {\it one-way feedback}, i.e., edges appear and disappear over time depending on the opinions of the two connected vertices, while the opinion dynamics is not affected by the graph structure. Our goal is to investigate the joint evolution of the entries of a {\it voter subgraph count vector}, i.e., vector of subgraphs where each vertex has a specific opinion, in the regime that the number of vertices grows large. The main result of this paper is a functional central limit theorem. In particular, we prove that, under a proper centering and scaling, the joint functional of the vector of subgraph counts converges to a specific multidimensional Gaussian process.

{\it Acknowledgment.} The work of NK was supported by the European Union’s Horizon 2020 research and innovation programme under the Marie Skłodowska-Curie grant agreement no.\ 101034253; SB received the support from the same grant while affiliated with Leiden University. SB was further supported through “Gruppo Nazionale per l’Analisi Matematica, la Probabilità e le loro Applicazioni” (GNAMPA-INdAM).
  \includegraphics[height=1em]{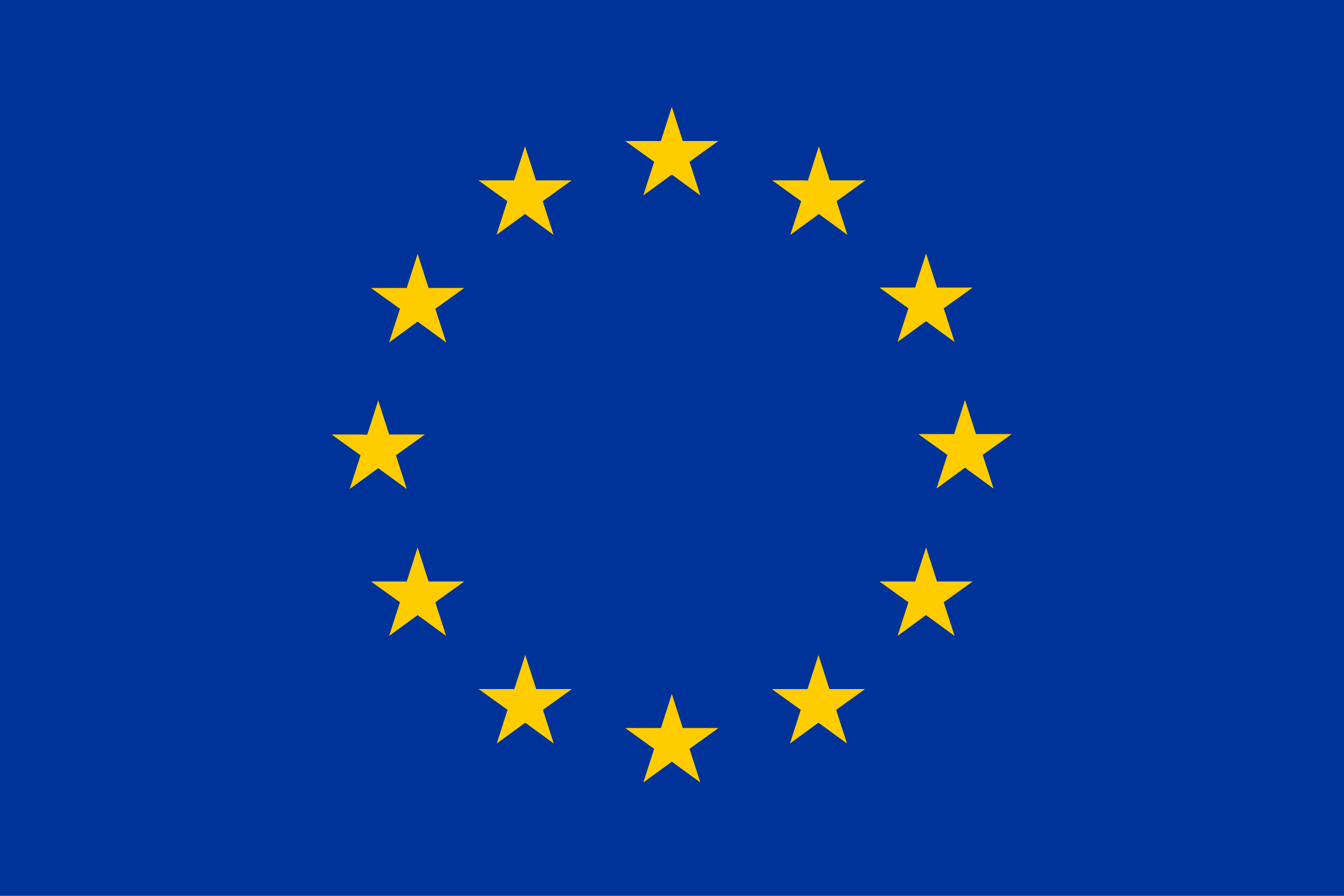}
\end{abstract}

\maketitle

\newpage

\section{Introduction}\label{Section:introduction}


The {\it voter model}, originally introduced as an interacting particle system on a lattice \cite{HL1975,L1985}, provides a fundamental framework for studying opinion dynamics in social networks. It describes how individuals update their opinions through local interactions, typically by adopting the state of a randomly chosen neighbor. Despite its simplicity, the model captures essential features of consensus formation, competition between different opinions, and the influence of network structure on collective behavior. 
Over the years, various extensions have been proposed to incorporate factors like external influences and heterogeneous interaction patterns, where individuals may have varying degrees of influence or interact at different rates depending on their position in the network. 

In recent years, the study of the voter model on {\it static} random graphs has become a key area of focus in network science (see e.g. \cite{CFR10,O13,ABHdHQ24,ACHQ24,C25}). Yet, most complex networks are not static but constantly {\it evolve}, with connections changing over time in response to external influences or internal interactions. This dynamic aspect introduces significant challenges in understanding the interplay between opinion formation and network evolution. Since such intertwined dynamics are prevalent in most real-world networks, developing a rigorous mathematical framework to study models that incorporate both elements is essential. However, despite significant efforts and empirical insights, our theoretical understanding remains limited. The few existing results are scattered across different settings, including studies on random walks on dynamic percolation \cite{PSS15,PSS20}, random walks on graphs with dynamic rewiring \cite{AGHdH18}, random walks on graphs with dynamic switching \cite{AGHdH19,ST20}, the contact process on dynamic random graphs \cite{JM17,CSW22}, the voter model on dynamic random graphs \cite{HN06,DGLMSSSV12,BS17,BBdHM2024coevolution,ABHdHQ25}. 

Functional Central Limit Theorems (FCLTs) are fundamental in probability theory as they provide insights into the fluctuations of graph statistics (e.g. subgraph counts, principal eigenvalue, degree sequence) by establishing their convergence to Gaussian processes. Understanding subgraph counts in static networks is crucial because they help characterize the graph's local and global structure, such as clustering or connectivity. They are used to compare real-world networks with random graph models, test hypotheses about network properties, and understand the asymptotic behavior of large graphs. While the classical theory of subgraph counts has been extensively studied for static graphs (see e.g. \cite{ER60,R88,JN91,O17,CDB24,BKL24}), introducing a dynamic component leads to additional stochastic dependencies that require new techniques. Recent progress in the study of FCLTs for dynamic random graphs has provided a solid foundation for our work. For instance, in \cite{RS16,BB2022PAM} the authors prove a FCLT for degree counts in preferential attachment models, while the focus of the authors in \cite{HKM2024,HKM2025FCLT} is the study of fluctuations for the principal eigenvalue and subgraph counts in a dynamic Erdős–Rényi graph.

In this paper, we focus on a two-state voter model on a dynamic random graph with one-way feedback between the joint dynamics of opinions and edges, namely, the network structure evolves based on opinions, while the dynamics of opinions remain unaffected by the graph. Indeed, in real-world scenarios, social ties are often influenced by individual opinions: people form and sever connections based on shared beliefs or ideological alignment, but, in many settings, these evolving connections do not immediately influence individual opinions. This assumption, while simplifying the analysis, still captures essential real-world scenarios in which individuals' opinions shape their social connections over time, but the structure of the network does not instantaneously influence opinion changes. 

A key novelty of our paper is that we extend the existing framework, recently developed in the context of FCLTs on dynamic random graphs, to incorporate opinion dynamics, thereby adding an additional layer of complexity due to the interaction between the evolving states of agents and the structural changes in the network. In particular, we establish a FCLT for the joint vector of subgraph counts of a voter model on dynamic random graphs, gaining deeper insights into the interplay between opinion dynamics and network structure, going beyond deterministic asymptotics to capture the inherent randomness of the system. We thus capture, for instance, the joint distribution of the number of edges between vertices agreeing on different opinions. To the best of our knowledge, this is the first paper to address this. 

Our work connects to the recent results in \cite{BBdHM2024coevolution}, where the convergence of dynamic random graphs to limiting graphon-valued processes is established. The connection between subgraph counts and graphon theory is well known: subgraph densities can indeed be derived from the graphon representation. In the graphon setting, there are several recent works investigating central limit theorems in the static case \cite{BCJ23,HPS21,KR21,N1989}, but, importantly, there seem to be no asymptotic normality results covering the setting of graphon-valued processes. Thus, here we describe how subgraph counts fluctuate around their expected values when the system evolves, which allows us to quantify the impact of network evolution on the variability of them. This is the main motivation of the paper.

To establish this FCLT, our proof combines several layers of analysis. We first identify the relevant fluctuation scale through careful moment computations, where the first and second moments reflect both the stochastic dynamics of the voter model and the temporal correlations induced by the evolving network. The convergence of finite-dimensional distributions then requires disentangling the interaction between opinion dynamics and edge dynamics, ensuring that the limiting covariance structure captures contributions from both sources of randomness. Finally, tightness is established by controlling short-time fluctuations of the subgraph counts, which relies on quantitative bounds that balance the speed of opinion updates with the rate of network rewiring. Taken together, these steps provide a framework that not only yields the desired functional convergence but also clarifies how the joint randomness of states and network evolution shapes the Gaussian limit.

The innovative contributions of the present paper can be summarized as follows.
\begin{enumerate}
	\item We establish a FCLT for the joint vector of subgraph counts in a voter model evolving on a dynamic random graph with one-way feedback.
	\item We rigorously characterize the limiting Gaussian process governing the fluctuations of subgraph counts in terms of its covariance function.
	\item Our analysis refines the understanding of subgraph count fluctuations going beyond deterministic asymptotics derived in \cite{BBdHM2024coevolution} in terms of convergence in the space of graphons.
\end{enumerate}

Our result shows that understanding the stochastic fluctuations of subgraph counts is indispensable for capturing how opinion dynamics and network evolution interact, reaching beyond deterministic approximations. This perspective is crucial in both theory and applications: in real systems, it is precisely these fluctuations, rather than averages alone, that can trigger consensus, sustain diversity, or drive abrupt shifts in collective behavior.

In the context of applications, a natural question is whether we can construct confidence bands for subgraph densities from the observed sample paths. This problem is closely related to statistical estimation in interacting particle systems on dynamic random graphs, for instance of the vertex and edge switching rates. In \cite{BW2025}, the authors address this question in the setting of a discrete-time version of the model considered here. They employ the method of moments to consistently estimate the parameters of the dynamics from partial information of the process, specifically from snapshots of the total number of edges present. Once estimates of the switching rates are available, they can in turn be used to estimate the limiting covariance function, which is a key ingredient in constructing confidence bands for subgraph densities. Extending these ideas to the continuous-time dynamics studied in the present work would provide a systematic statistical framework for connecting observed sample paths with the asymptotic fluctuation theory. This extension would play a central role in applications, by linking theoretical fluctuation results with practical tools for uncertainty quantification and calibration.

Although this paper focuses on a one-way feedback interaction between opinion and edge dynamics, our study provides a fundamental step toward understanding the rich interplay between opinion dynamics and structural evolution in networks. A natural and significant extension would be indeed to consider fully {\it co-evolutionary} models, where opinions and network structure influence each other bidirectionally. At the end of this paper, we provide a heuristic argument explaining why the characterization of the limiting distribution of the subgraph count process is significantly more complex when dealing with a two-way feedback model. Investigating the effects of such mutual feedback on subgraph count fluctuations remains an important open problem for future research. 

The remainder of the paper is organized as follows. Section \ref{Section:one-way_feedback} is devoted to the presentation of the model of interest, the required notation, and the main contribution of this work. In particular, in Section \ref{sec:defmodel} we introduce the voter model on a dynamic random graph with one--way feedback. Section \ref{Section:graph_notations} contains some graph-related notations, while in Section \ref{Section:main_results} we state our main result concerning the functional central limit theorem for the joint subgraph count process in Theorem \ref{main_model_1}. Finally, in Section \ref{sec:discussion} we highlight how the main result of the present paper extends the main findings in \cite{BBdHM2024coevolution} for the same model. Section \ref{Section:proofs} contains the proof of Theorem \ref{main_model_1}. In particular, in Sections \ref{Section:proof_of_covariance}--\ref{Section:proof_of_tightness} we prove the three key lemmas that build up the proof of the main theorem, while in Section \ref{Section:proofprop} we provide a characterization of the variance function for the limiting Gaussian process. Finally, in Section \ref{Section:two-ways_feedback} we introduce a voter model on a dynamic random graph with a two-way feedback interaction between the joint dynamics of vertices and edges, and we provide some arguments explaining why the limiting distribution of the subgraph count process is much more challenging to derive. A complete characterization is left as a direction for future research.

\section{One-way feedback model: definition and main result}\label{Section:one-way_feedback}

Throughout the paper we use the following notations. For any natural number $n$, we write $\scale{n}$ to denote the set $\{1,2,...,n\}$. All vectors are written in bold and, vice versa, all bold symbols refer to vectors. All entries of a vector share the same letter, but do not use bold font and have a subscript which represents the index of the entry. For example, if $\vk b\in\R^d$, then $\vk b=(b_1,...,b_d)^{\top}$. If a vector has a subscript itself, for each element of this vector we first write a subscript of the vector, and after the subscript of the element. For example, if $\vk b_i\in\R^d$, then $\vk b_i=(b_{i,1},...,b_{i,d})^{\top}$. Finally, for any finite set $A$, we denote the cardinality of $A$ as $|A|$.

\subsection{Definition of the model}\label{sec:defmodel}

We consider a model with a one-way interaction (i.e., the states of the vertices are independent of the renewal rates of the edges, but the reverse is not true, meaning that the states of the edges depends on the renewal rates of the vertices). For a given deterministic time horizon $T$, let $(G_n(t))_{t\in[0,T]}$ denote our graph-valued process. We restrict ourselves to consider simple graphs, i.e., graphs without self-loops and multiple edges. At any time $t$, the set of vertices is $V=\scale{n}$ and any vertex holds opinion $+$ or $-$. We assume that each edge is initially active with probability $p_0$, independently of everything else (i.e., we start from an Erdös-Rényi random graph on $n$ vertices with ‘edge retention probability' $p_0$). The model is characterised by the following dynamics.

\smallskip
{\it Vertex dynamics.} Each vertex is re-sampled at a rate that depends on its opinion, i.e., a rate-$\gamma_{-+}$ (resp.\ $\gamma_{+-}$) Poisson clock is attached to each vertex holding opinion $-$ (resp.\ $+$). When the clock rings, the associated vertex switches its opinion from $-$ to $+$ (resp.\ from $+$ to $-$). All happens independently of everything else.

\smallskip
{\it Edges dynamics.} Each edge is re-sampled at rate $1$, i.e., a rate-$1$ Poisson clock is attached to each edge and when the clock rings the edge is active with a probability that depends on the current opinion of the two connected vertices: with probability $\pi_+$ if the two vertices hold opinion $+$, with probability $\pi_-$ if the two vertices hold opinion $-$, and otherwise with probability $\tfrac12(\pi_++\pi_-)$.

\smallskip
By definition, it is immediate to see that the vertex dynamics is self-sufficient, meaning that it is independent of the Poisson clocks attached to the edges. However, it is not true for the edge dynamics, which depends on the opinions of the vertices through $\pi_+$ and $\pi_-$. Throughout this paper, we associate each active (resp.\ inactive) edge with 1 (resp.\ 0), so that the resulting edge process for the vertex-pair $(u,v)$, which we denote by $(a_{u,v}(t))_{t\in[0,T]}$, is
\begin{equation}\label{eq:edgeprocess}
	a_{u,v}(t) =
	\begin{cases}
		1 &\text{ if } (u,v) \text{ is active in } G_n(t), \\
		0 &\text{ if } (u,v) \text{ is inactive in } G_n(t).
	\end{cases}
\end{equation}

\subsection{Graph notations}\label{Section:graph_notations}

In this section we introduce various graph-related notations.

\begin{itemize}

	\item We define a {\it voter graph} as a graph in which each vertex is associated with an opinion, which can be $+$ or $-$. We denote the opinion of vertex $v$ of a voter graph $H$ as $x_v(H)\in\{+,-\}$. Moreover, for any voter graphs $H$ and $H^{\prime}$ we call them isomorphic only if they are isomorphic as simple graphs and this isomorphism preserves the opinions of the vertices in each graph. 
	
	\item For any voter graph $H$, let $\mathcal{V}(H)$, $\mathcal{E}(H)$ and $\mathcal{A}(H)$ be the number of vertices, edges and automorphisms (which preserves the opinions) of the graph $H$, respectively. 
	
	\item We say that a graph $H$ is {\it labeled} if its vertices have different integer labels. We denote by ${\tt V}(H)\subset \N$ the set of labels of vertices of a labeled graph $H$, and by ${\tt E}(H)\subset\{(u,v)\in\N^2\colon u<v\}$ the set of its edges. Two labeled graphs $H,H^{\prime}$ are considered the same if ${\tt V}(H) = {\tt V}(H^{\prime})$ and ${\tt E}(H) = {\tt E}(H^{\prime})$. 
	
	\item For any ${\tt S}\subset\N$, we define ${\tt LVG}_{{\tt S}}$ as the set of all labeled voter graphs $h$ such that ${\tt V}(h)\subset{\tt S}$. Throughout the paper we assume, for ${\tt S}_1\subset{\tt S}_2\subset\N$, the natural inclusion 
	\begin{align*}
		{\tt LVG}_{{\tt S}_1}\subset{\tt LVG}_{{\tt S}_2}.
	\end{align*}
	
	\item For any ${\tt S}\subset\N$ and any voter graph $H$, we define the set ${\tt G}_{{\tt S}}(H)$ of all different labeled voter graphs from ${\tt LVG}_{{\tt S}}$ which are isomorphic to $H$, i.e.,
	\begin{align}
		{\tt G}_{{\tt S}}(H) = \{h\in {\tt LVG}_{{\tt S}}\colon h \text{ is isomorphic to } H \}.\label{tt_G_def}
	\end{align}
	\begin{remark}
		Recall that the labels of the graph $h$ are ignored in the isomorphism used in the definition of ${\tt G}_{{\tt S}}(H)$.
	\end{remark}
	\item For any $h,h^{\prime}\in{\tt LVG}_{{\scale{n}}}$, we write $h\sqsubset h^{\prime}$ if ${\tt V}(h)\subset {\tt V}(h^{\prime})$ and ${\tt E}(h)\subset {\tt E}(h^{\prime})$.
\end{itemize}

Given a voter graph $H$, for any set $S\subset\N$ such that $\abs{{\tt S}}\geqslant \mathcal{V}(H)$, we get
\begin{align}
	\abs{{\tt G}_{{\tt S}}(H)} = \frac{\abs{{\tt S}}!}{(\abs{{\tt S}}-\mathcal{V}(H))!\mathcal{A}(H)}.\label{G_cardinality}
\end{align}

Note that $G_n(t)$ is a random element of ${\tt LVG}_{\scale{n}}$ for any $t\in[0,T]$. So, for any $h\in {\tt LVG}_{{\scale{n}}}$ we can define the indicator function which tells us if  $h$ is a labeled subgraph of $G_n(t)$ as
\begin{align}
	\mathbb{I}_n(h,t) = \mathbb{I}\{h \sqsubset G_n(t)\} = \prod_{(u,v)\in{\tt E}(h)}a_{u,v}(t)\prod_{v\in{\tt V}(h)}\mathbb{I}\{x_{v}(G_n(t)) = x_{v}(h)\}, \label{I_N_def}
\end{align}
where we recall that $a_{u,v}(t)$ is defined in \eqref{eq:edgeprocess}. Using the symmetry of the model with respect to any permutation of the vertices, we obtain that the distribution of the random variable $\mathbb{I}_n(h,t)$ does not depend on the particular choice of $h\in{\tt G}_{\scale n}(H)$ for any voter graph $H$.

For any ${\tt S}\subset\N$ and $k\geq2$, we define the space of vectors of length $k$ with ordered elements in ${\tt S}$ as
\begin{align*}
	{\tt K}^k({\tt S}) = \{\vk v\in {\tt S}^{k}\colon v_1<\ldots<v_k\}.
\end{align*}
Throughout the paper we use a natural bijection between vectors $\vk v\in{\tt K}^k({\tt S})$ and subsets of ${\tt S}$ of cardinality $k$. So, we can treat every vector $\vk v$ simultaneously as a vector in ${\tt K}^k({\tt S})$ and as a subset of ${\tt S}$. 

For further calculations, the constants defined below will play a crucial role. For any voter graph $H$ and any $h\in{\tt G}_{\scale{n}}(H)$, we define
\begin{align}
	\mathcal{P}_{H}(t) &= \E{\mathbb{I}_{n}(h,t)} \label{P_def},
\end{align}
and for any two voter graphs $H,H^{\prime}$ and any $\vk v_1\in{\tt K}^{\mathcal{V}(H)}(\scale{n})$, $\vk v_2\in{\tt K}^{\mathcal{V}(H^{\prime})}(\scale{n})$ such that $\abs{\vk v_1\cap \vk v_2} = 1$, we define
\begin{align}
	\mathcal{C}_{H,H^{\prime}}(s,t) &= \frac{1}{\bigl((\mathcal{V}(H)-1)!\bigr)\bigl((\mathcal{V}(H^{\prime})-1)!\bigr)}\sum_{\substack{h\in {\tt G}_{\vk v_1}(H) \\ h^{\prime}\in {\tt G}_{\vk v_2}(H^{\prime})}}\operatorname{Cov}\bigl(\mathbb{I}_{n}(h,s),\,\mathbb{I}_{n}(h^{\prime},t)\bigl).\label{C_def}
\end{align}
Recall that using symmetry of our model, the constant $\mathcal{P}_H(t)$ does not depend on $n$ and the choice of the subgraph $h$, and the constant $\mathcal{C}_{H,H^{\prime}}(t,s)$ does not depend on $n$ and the choice of vectors $\vk v_1$ and $\vk v_2$. For simplicity, later we will use $\mathcal{C}_{H}(s,t)$ to denote $\mathcal{C}_{H,H}(s,t)$.

\subsection{Main result}\label{Section:main_results}

In this section we state the main result of this paper, which is the functional central limit theorem for the simultaneous subgraph count process. To this end, for a given $m\in\N$, we fix $m$ deterministic voter graphs $H_1,\ldots,H_m$. For any $n\in\N$, $i\in\scale{m}$ and $t\in[0,T]$, we define $X_{n,i}(t)$ as the number of different inclusions of the graph $H_i$ in $G_n(t)$, which preserve the opinions of each vertex in the the graph $H_i$. Alternatively, the process $X_{n,i}(\cdot)$ can be defined, for $t\in[0,T]$, as
\begin{align}
	X_{n,i}(t)=\sum_{h\in{\tt G}_{\scale{n}}(H_i)}\mathbb{I}_n(h,t),\label{X_def}
\end{align}
where ${\tt G}_{\scale{n}}(H_i)$ and $\mathbb{I}_n(h,t)$ are defined in \eqref{tt_G_def} and \eqref{I_N_def}, respectively. The main objective of this paper is to analyze the asymptotic behavior of the following $m$-dimensional stochastic process
\begin{align}
	\vk X_n(\cdot) = \bigl(X_{n,1}(\cdot),\ldots, X_{n,m}(\cdot)\bigr)^{\top}. \label{vk_X_def}
\end{align}We state below the main result of this contribution.
\begin{theo}\label{main_model_1}
	Let $G_n(t)$ for $t\in[0,T]$ be a one-way feedback voter model on a dynamic random graph with $n$ vertices (as described in Section \ref{sec:defmodel}), and let $H_1,\ldots,H_m$ be deterministic voter graphs for some $m\in\N$. Then, as $n\to\infty$, the process $\vk X_n(\cdot)$ defined in \eqref{vk_X_def} satisfies the following limiting theorem
	\begin{align}
		\left(\frac{X_{n,1}(t)-\E{X_{n,1}(t)}}{n^{\mathcal{V}(H_1)-1/2}},\ldots,\frac{X_{n,m}(t)-\E{X_{n,m}(t)}}{n^{\mathcal{V}(H_m)-1/2}}\right)\to \vk X(t), \label{eq:main_theo}
	\end{align}
	in distribution in $\mathbb{D}([0,T],\R^m)$, where $\vk X(\cdot)$ is a centered vector-valued Gaussian process characterized, for any $i,j\in\scale{m}$ and any $0\leq s \leq t \leq T$, by the covariance matrix function $\Sigma(\cdot,\cdot)$, where for any $s,t\in[0,T]$, $i,j\in\scale{m}$
	\begin{align*}
		\Sigma(s,t)_{i,j} = \mathcal{C}_{H_i,H_j}(s,t),
	\end{align*}
	where the constant $\mathcal{C}_{H_i,H_j}(s,t)$ is defined in \eqref{C_def}.
\end{theo}

{\begin{remark}
    For any $n\in\N$, $i\in\scale{m}$, $t\in[0,T]$, the exact formula for $\E{X_{n,i}(t)}$ is provided in \eqref{exp_formula}, in combination with \neprop{prop:variance}, while the exact formula for $\Sigma(s,t)$ for any $s,t\in[0,T]$ is provided in \eqref{eq:C_G(t,t)}. From \eqref{cov_formula} we can see that $\Sigma(\cdot,\cdot)$ is non-negative definite (as a limit of non-negative definite covariance matrix functions), and \neprop{prop:variance} shows that $\Sigma(t,t)_{i,i}>0$ for any $i\in\scale{m}$, $t\in[0,T]$. Hence, $\Sigma(\cdot,\cdot)$ is a valid non-degenerate covariance matrix function.
\end{remark}

{\begin{remark}\label{rmk:centering} Since $\E{X_{n,i}(t)}$ can be calculated explicitly for any $n\in\N$, $i\in\scale{m}$, and $t\in [0,T]$, and satisfy the asymptotic relation as $n\to\infty$ uniformly in $t$ (see \nelem{exp_and_cov} below for more details)
\begin{align*}
    \E{X_{n,i}(t)} = n^{\mathcal{V}(H_i)}\frac{\mathcal{P}_{H_i}(t)}{\mathcal{A}(H_i)} + o\left(n^{\mathcal{V}(H)-1/2}\right),
\end{align*}
the limiting theorem \eqref{eq:main_theo} can be equivalently written as
\begin{align*}
		\left(\frac{X_{n,1}(t)-n^{\mathcal{V}(H_1)}\frac{\mathcal{P}_{H_1}(t)}{\mathcal{A}(H_1)}}{n^{\mathcal{V}(H_1)-1/2}},\ldots,\frac{X_{n,m}(t)-n^{\mathcal{V}(H_m)}\frac{\mathcal{P}_{H_m}(t)}{\mathcal{A}(H_m)}}{n^{\mathcal{V}(H_m)-1/2}}\right)\to \vk X(t).
	\end{align*}
\end{remark}

\subsection{Discussion}\label{sec:discussion}

The model defined in Section \ref{sec:defmodel} was introduced in \cite{BBdHM2024coevolution}, where the authors generalized the results in \cite{BdHM2022} by proving a functional law of large numbers for a voter model on a dynamic random graph in the space of graphons, namely, the empirical graphon dynamics associated with the graph-valued process $(G_n(t))_{t\in[0,T]}$ converges to a deterministic limit process described by a graphon process. (See \eqref{eq:empgraphon} for the definition of empirical graphon associated with a graph). In this section, we highlight the relationship between the main findings in \cite{BBdHM2024coevolution} and Theorem \ref{main_model_1}. Furthermore, we show how the result in the present paper refine and extend those in \cite{BBdHM2024coevolution}, providing a more detailed and comprehensive analysis.

Any finite simple undirected graph $G_n$ on $n$ vertices can be represented as a graphon $\mathfrak{h}^{G_n}$ by setting
\begin{equation}\label{eq:empgraphon}
	\mathfrak{h}^{G_n}(x,y) = 
	\begin{cases}
		1 &\text{ if there is an edge between vertex } \lceil nx \rceil \text{ and vertex } \lceil ny \rceil, \\
		0 &\text{ otherwise},
	\end{cases}
\end{equation}
which we refer to as {\it empirical graphon} associated with $G_n$. Given a finite simple graph $H$, the subgraph density of $H$ in $G_n$ is given by
\[
t(H, G_n) = \int_{[0,1]^{\mathcal{V}(H)}} \prod_{(u,v)\in {\tt E}(H)} \mathfrak{h}^{G_n}(x_u,x_v) \prod_{u \in {\tt V}(H)} \dd x_u.
\]
\cite[Theorem 2.2]{BBdHM2024coevolution} shows that, as $n\to\infty$, $(\mathfrak{h}^{G_n(t)})_{t\in[0,T]}$ converges to the deterministic graphon process $(\mathfrak{g}^{[F]}(t)_{t\in[0,T]}$, which can be described in terms of the function 
\begin{equation}
	\label{eqn:Hdef}
	\mathcal{H}(t;u,v) = p_0 {\rm e}^{-t}+ \dfrac{1}{2}\left[ \pi_+  u + \pi_- (1-{\rm e}^{-t} - u) 
	+  \pi_+  v +  \pi_- (1-{\rm e}^{-t} - v) \right],
\end{equation}
where $t\in[0,T]$ and $u,v\in[0,1]$. See \cite[eq.\ (1.11)]{BBdHM2024coevolution} for the precise definition of $\mathfrak{g}^{[F]}$. Thus, as $n\to\infty$, the subgraph density of $H$ in $G_n$ at time $t^*$ is given by
\[
t(H, \mathfrak{g}^{[F]}(t^*)) = \int_{[0,1]^{\mathcal{V}(H)}} \prod_{(u,v)\in {\tt E(H)}} \mathfrak{g}^{[F]}(t^*;x_u,x_v) \prod_{u \in {\tt V}(H)} \dd x_u.
\]
This result enables us to express, as $n\to\infty$, the total number $X_{n,H}(t^*)$ of subgraphs of type $H$ in $G_n$ at time $t^*$ in terms of the graphon-based density as
\[
\lim_{n\to\infty} X_{n,H}(t^*) = \lim_{n\to\infty} t(H, \mathfrak{g}^{[F]}(t^*)) n^{\mathcal{V}(H)}.
\]
Note that this term, being the leading one appearing in the LLN, coincides with the explicit centering in Remark \ref{rmk:centering}, so that it is possible to replace the centering in the FCLT by the limiting subgraph density of the graphon process. This relationship illustrates how subgraph counts can be derived from graphon-based densities, providing a crucial link between large-scale graph structure and subgraph statistics.

The main contribution of this paper is then the characterization of the fluctuations around this limiting value. Specifically, Theorem \ref{main_model_1} shows that the centered and scaled joint vector of simultaneous subgraph count converges to a multidimensional centered Gaussian process, thereby providing a functional CLT that describes the probabilistic nature of these fluctuations. The FCLT derived here strengthens the graphon-based asymptotic analysis, emphasizing its role in characterizing the statistical properties of large random graphs. The covariance function of the limiting Gaussian process encodes crucial structural information about the interdependence between different subgraph counts, which is fundamental to understand correlations in large graph structures. In what follows we recall some quantities introduced in \cite{BBdHM2024coevolution} in order to provide an explicit formula for the covariance function. Due to the intricate nature of the double dynamics of vertices and edges, this formula cannot be written in closed form. In Proposition \ref{prop:variance} we provide such a representation for the variance function. Similar arguments can be applied to derive the covariance function as well, but since the computations are long and technical, and mostly not inspiring in terms of tools that we need, we do not push further on this point. Even if this formula has no closed form, its explicit formulation allows us to make quantitative predictions about the behavior of subgraph counts, offering insights into the intricate dependence on the parameters of the model.

A key quantity we need in what follows is the {\it type} of vertex $i$ at time $t$, which encodes how much time vertex $i$ has held opinion $+$ up to time $t$, and it is formally defined as 
\begin{equation}\label{eqn:typedef}
y_i(t) = \int_0^t \dd s \eee^{-s} \mathbb{I}\{x_i(G_n(t-s))=+\}.
\end{equation}
This {\it specific} functional allows us to express the probability that edge $ij$ is active at time $t$ in terms of $y_i(t)$ and $y_j(t)$ only. In these terms, the function $\mathcal{H}(t;u,v)$ defined in \eqref{eqn:Hdef} represents indeed the probability that there is an active edge connecting two vertices with type $u$ and $v$ at time $t$ (see \cite[Section 1.2]{BBdHM2024coevolution} for more details on how the vertex type should be carefully chosen). This means that keeping track of the opinions of vertices {\it only} is not enough to capture the structure of the resulting graph. We then need to introduce the {\it generalised type} of vertex $i$ at time $t$ as the pair $X_i(G_n(t)) = (x_i(G_n(t)), y_i(t))$, which belongs to the space $\{-,+\}\times [0,1]$, and whose density is referred to as
\[
f_+(t,u) \dd u = \mathbb{P}(X_i(G_n(t))\in(+,\dd u)), \qquad f_-(t,u) \dd u = \mathbb{P}(X_i(G_n(t))\in(-,\dd u)).
\]
The densities $f_+(\cdot,\cdot)$ and $f_-(\cdot,\cdot)$ satisfy a system of PDEs describing their time evolution, which allows for a characterization of these densities, albeit an implicit one (see \cite[Theorem 2.3]{BBdHM2024coevolution} for more details).

We are now ready to provide the explicit characterization of the variance function of the limiting Gaussian process, which is the content of the next proposition.
While the expression in equation \eqref{eq:C_G(t,t)} below is technically involved, its underlying idea can be summarized as follows. At any time $t$, any edge $ij$ is conditionally independent of all the other edges given the types of the vertices, and the probability of being active can be indeed expressed in terms of the types of vertices $i$ and $j$ {\it only} (see \eqref{eqn:Hdef}). This means that, when computing the edge probabilities at time $t$, it is possible to factorize them given the types of all the involved vertices. This is the key tool we use to derive the formulas appearing in the following proposition, whose proof is deferred to Section \ref{Section:proofprop}.

\begin{prop}\label{prop:variance}
For any $t\in[0,T]$, any two voter graphs $H,H^{\prime}$ and any $\vk v_1\in{\tt K}^{\mathcal{V}(H)}(\scale{n})$, $\vk v_2\in{\tt K}^{\mathcal{V}(H^{\prime})}(\scale{n})$ such that $\vk v_1\cap \vk v_2 = \{k\}$, 
   	 \begin{align}\label{eq:C_G(t,t)}
			\mathcal{C}_{H, H^{\prime}}(t,t) &= \displaystyle \frac{1}{\bigl((\mathcal{V}(H)-1)!\bigr)\bigl((\mathcal{V}(H^{\prime})-1)!\bigr)} \sum_{\substack{h\in {\tt G}_{\vk v_1}(H) \\ h^{\prime}\in {\tt G}_{\vk v_2}(H^{\prime})}} \Biggl\{ \int_0^1 \Biggl[ \displaystyle \int_{[0,1]^{|h|-1}} \prod_{(u,v)\in{\tt E}(h)} \mathcal{H}(t;\ell_u,\ell_v) \notag \\ 
			&\quad \times \displaystyle \prod_{i\in {\tt V}(h)\setminus\{k\}} f(t,\ell_i) \dd \ell_i \int_{[0,1]^{|h^{\prime}|-1}} \prod_{(w,z)\in{\tt E}(h^{\prime})} \mathcal{H}(t;\ell_w,\ell_z) \notag \\
			&\quad \times \displaystyle \prod_{j\in {\tt V}(h^{\prime})\setminus\{k\}} f(t,\ell_j) \dd \ell_j \sum_{o\in\{-,+\}} \Biggl( \prod_{i\in{\tt V}(h)} \mathbb{P}(x_{i}(G_n(t)) = x_{i}(h)| x_k(G_n(t))=o) \notag \\
			&\quad \displaystyle \times \prod_{j\in{\tt V}(h^{\prime})} \mathbb{P}(x_{j}(G_n(t))= x_{j}(h^{\prime})| x_k(G_n(t))=o) \int_0^1 f_o(t,\ell_k) \dd \ell_k \Biggr) \Biggr] f(t,\ell_k) \dd \ell_k \notag \\
			&\quad - \mathcal{P}_{H}(t) \mathcal{P}_{H'}(t) \Biggr\},
	\end{align}
    where $f(t,\ell)=f_+(t,\ell)+f_-(t,\ell)$ for any $\ell\in[0,1]$ and, for any voter graph $\bar{H}$ and any $\bar h\in{\tt G}_{\scale{n}}(\bar{H})$,
    \[
    \mathcal{P}_{\bar{H}}(t) = \displaystyle \int_{[0,1]^{|\bar{h}|}} \prod_{(u,v)\in{\tt E}(\bar{h})} \mathcal{H}(t;\ell_u,\ell_v) \prod_{i\in {\tt V}(\bar{h})}  \mathbb{P}(x_{i}(G_n(t)) = x_{i}(\bar{h})) f(t,\ell_i) \dd \ell_i.
    \]
	In particular, for any $t\in[0,T]$, $\mathcal{C}_H(t,t)>0$.
\end{prop}

\section{Proof of the main result}\label{Section:proofs}

\prooftheo{main_model_1}
To prove \netheo{main_model_1}, we follow the three steps below.
\begin{description}
	\item[Step 1] Computation of expectation and covariance functions of $\vk X_n(\cdot)$.
	\item[Step 2] Establish finite-dimensional convergence of distributions.
	\item[Step 3] Show tightness that yields the desired functional convergence.
\end{description}

The first step is based on the following lemma, whose proof is deferred to Section \ref{Section:proof_of_covariance}.
\begin{lem}\label{exp_and_cov}
	For any $n\in\N$, $i\in\scale{m}$ and $t\in[0,T]$,
	\begin{align}
		\E{X_{n,i}(t)} = \frac{n! \mathcal{P}_{H_i}(t)}{(n-\mathcal{V}(H_i))!\mathcal{A}(H_i)},\label{exp_formula}
	\end{align}
	and for any $i,j\in\scale{m}$,  $s,t\in[0,T]$, as $n\to\infty$,
	\begin{align}
		\covv{X_{n,i}(s),\,X_{n,j}(t)} = n^{\mathcal{V}(H_i) + \mathcal{V}(H_j) - 1}\left(\mathcal{C}_{H_i,H_j}(s,t)+o(1)\right).\label{cov_formula}
	\end{align}
\end{lem}

The second step consists in showing the convergence of finite dimensional distributions.  For this purpose, it is significantly more convenient to work with the centered and normalized version $\vk X^{\star}_n(\cdot)$ of the initial process $\vk X_n(\cdot)$, which can be defined for any $i\in\scale{m}$, $t\in[0,T]$ as
\begin{align}
	X^{\star}_{n,i}(t) = \frac{X_{n,i}(t)-\E{X_{n,i}(t)}}{n^{\mathcal{V}(H_i)-1/2}}.\label{X_star_def}
\end{align}
Note that the process $\vk X^{\star}_n(\cdot)$ is completely characterized thanks to \nelem{exp_and_cov}. Then, convergence of finite dimensional distributions of $\vk X_n(\cdot)$ is equivalent to the following lemma, whose proof is postponed to Section \ref{Section:proof_fdd_convergence}.
\begin{lem}\label{asymptotic_normality}
	For any $k\in\N$ and $t_1,\ldots,t_k\in[0,T]$, the random vector $\bigl(\vk X^{\star}_{n}(t_1),\ldots, \vk X^{\star}_{n}(t_{k})\bigr)$ converges in distribution to a centered Gaussian random vector with covariance block matrix $\vk \Sigma$ with elements $\left(\Sigma_{u,v}\right)_{i,j} = \mathcal{C}_{H_i,H_j}(t_u,t_v)$ for $i,j\in \scale{m}$, $u,v\in\scale{k}$, where the constant $\mathcal{C}_{H_i,H_j}(t_u,t_v)$ is defined in \eqref{C_def}.
\end{lem}

Finally, to justify the claim of \netheo{main_model_1}, we apply the multidimensional analogue of \cite[Theorem 13.5]{billingsley2013convergence}, in particular, of the condition given in \cite[eq.\ (13.14)]{billingsley2013convergence}. According to this condition, to verify that the sequence  $\vk X_n^{\star}(\cdot)$ is tight in $\mathbb{D}([0,T],\R^m)$ it is sufficient to find a continuous function $F(\cdot)$ such that, for any $0\leqslant r<s<t\leqslant T$,
\begin{align}
	\E{\abs{\vk X^{\star}_n(t)-\vk X^{\star}_n(s)}^2\abs{\vk X^{\star}_n(s)- \vk X^{\star}_n(r)}^2}\leq \abs{F(t)-F(r)}^2.
\end{align}
The required function $F(\cdot)$ is characterized in the following lemma, whose proof is deferred to Section \ref{Section:proof_of_tightness}.

\begin{lem}\label{lem:tightness} For any $n\in\N$ and any $0\leqslant r<s<t\leqslant T$,
	\begin{align}
		\E{\abs{\vk X^{\star}_n(t)- \vk X^{\star}_n(s)}^2\abs{\vk X^{\star}_n(s) - \vk X^{\star}_n(r)}^2}\leq \left(\sum_{i,j=1}^{m} \left(\bigl(2\mathcal{V}_{i,j}-2\bigr)!\right)^3F^2_{i,j}\right)\abs{t-r}^2,\label{tightness_claim}
	\end{align}
	where for any $i,j\in\scale{m}$
	\begin{align}
		F_{i,j} &= 4\max(\gamma_{-+},\gamma_{+-},\pi_+,\pi_-)\bigl(\mathcal{V}(G_i) + \mathcal{E}(H_i)+ \mathcal{V}(H_j)+ \mathcal{E}(H_j)\bigr),\label{F_ij_def}\\
		\mathcal{V}_{i,j}&= \mathcal{V}(H_i) + \mathcal{V}(H_j). \notag
	\end{align}
\end{lem}
Finally, the convergence claimed in \netheo{main_model_1} can be obtained by combining \cite[Theorem 13.5]{billingsley2013convergence} with Lemmas \ref{asymptotic_normality} and \ref{lem:tightness}.
\QED

\subsection{Proof of \nelem{exp_and_cov}}\label{Section:proof_of_covariance}

The goal of this section is to prove \eqref{exp_formula} and \eqref{cov_formula}. We focus first on \eqref{exp_formula}. Taking expectation of both sides in \eqref{X_def}, for any $i\in\scale{m}$,  $n>\mathcal{V}(H_i)$ and $t\in[0,T]$, we obtain
\[
\begin{array}{ll}
	\E{X_{n,i}(t)} &= \displaystyle\E{\sum_{h\in{\tt G}_{\scale{n}}(H_i)}\mathbb{I}_n(h,t)} =\sum_{h\in{\tt G}_{\scale{n}}(H_i)}\E{\mathbb{I}_n(h,t)}\\
	&= \displaystyle \sum_{h\in{\tt G}_{\scale{n}}(H_i)}\mathcal{P}_{H_i}(t)=\frac{n!}{(n-\mathcal{V}(H_i))!\mathcal{A}(H_i)}  \mathcal{P}_{H_i}(t),
\end{array}
\]
where in the last line we first used \eqref{P_def}, and then \eqref{G_cardinality}. Hence, \eqref{exp_formula} follows.

We proceed now to prove \eqref{cov_formula}. We set $V^{\star}_{i,j}=\min(\mathcal{V}(H_i),\mathcal{V}(H_j))$ and, for any $S\subset \N$, 
\begin{align*}
	{\tt G}_{S}^{i,j} &= {\tt G}_{S}(H_i)\times {\tt G}_{S}(H_j),\\
	{\tt G}_{S}^{i,j,k} &= \left\{\vk h\in{\tt G}_{S}^{i,j}\colon \abs{{\tt V}(h_1)\cap {\tt V}(h_2)}=k\right\}.
\end{align*}
Using again \eqref{X_def}, we can write the elements of the covariance matrix function of $\vk X_n(\cdot)$ for any $n\in\N$, any $i,j\in\scale{m}$ and any $s,t\in[0,T]$ as
\begin{align}
	\covv{X_{n,i}(s),\,X_{n,j}(t)} &= \Cov{\sum_{h\in{\tt G}_{\scale{n}}(H_i)}\mathbb{I}_n(h,s),\,\sum_{h\in{\tt G}_{\scale{n}}(H_j)}\mathbb{I}_n(h,t)} \notag \\
	&=\sum_{\vk h\in {\tt G}_{\scale{n}}^{i,j}}\covv{\mathbb{I}_n(h_1,s),\,\mathbb{I}_n(h_2,t)} \notag \\
	&=\sum_{k=0}^{V^{\star}_{i,j}}\sum_{\vk h\in {\tt G}_{\scale{n}}^{i,j,k}}\covv{\mathbb{I}_n(h_1,s),\,\mathbb{I}_n(h_2,t)}. \label{cov_repr_0}
\end{align}
We use \eqref{cov_repr_0} as follows: we show that only the term corresponding to $k=1$ in the inner sum contributes to the asymptotic of the covariance, and the rest terms are asymptotically negligible. Consider first the inner sum in the case $k=0$. Thus, ${\tt V}(h_1)\cap {\tt V}(h_2) = \varnothing$, and, by the definition of the model, the random variables $\mathbb{I}_n(h_1,s)$ and $\mathbb{I}_n(h_2,t)$ are independent. So, in \eqref{cov_repr_0} we can get rid of the sum for $k=0$, leading to
\begin{equation}
	\covv{X_n(s),\,X_n(t)} = \sum_{\vk h\in {\tt G}_{\scale{n}}^{i,j,1}}\covv{\mathbb{I}_n(h_1,s),\,\mathbb{I}_n(h_2,t)} \label{cov_repr_1}+ \sum_{k=2}^{V_{i,j}^{\star}}\sum_{\vk h\in {\tt G}_{\scale{n}}^{i,j,k}}\covv{\mathbb{I}_n(h_1,s),\,\mathbb{I}_n(h_2,t)}.
\end{equation}
Denote for simplicity 
\begin{align}
	V^k_{i,j} := \mathcal{V}(H_i) + \mathcal{V}(H_j) - k.\label{V_ij_def}
\end{align}
We calculate now the asymptotics of the first term of \eqref{cov_repr_1}, by fixing first the set ${\tt V}(h_1)\cup {\tt V}(h_2)$, then separating it into sets ${\tt V}(h_1)$ and ${\tt V}(h_2)$, and finally using \eqref{C_def} as follows (where for any $n\in\N$, $a,b\in\scale{n}$ and $u\in {\tt K}^{a}(\scale{n})$, $v\in {\tt K}^{b}(\scale{n})$ we denote by ${\tt I}(\vk u,\vk v)$ the number of common elements (disregarding their position) of vectors $\vk u$ and $\vk v$):
\begin{align}
	\sum_{\vk h\in {\tt G}_{\scale{n}}^{i,j,1}}\covv{\mathbb{I}_n(h_1,s),\,\mathbb{I}_n(h_2,t)}&= \sum_{\vk v\in{\tt K}^{V_{i,j}^1}(\scale{n})}\sum_{\vk h\in {\tt G}_{\vk v}^{i,j,1}}\covv{\mathbb{I}_n(h_1,s),\,\mathbb{I}_n(h_2,t)\}}\notag\\
	&= \sum_{\vk v\in{\tt K}^{V_{i,j}^1}(\scale{n})}\sum_{\substack{\vk{v}_1\in{\tt K}^{\mathcal{V}(H_i)}(\vk{v})\\\vk{v}_2\in{\tt K}^{\mathcal{V}(H_j)}(\vk{v}) \\ {\tt I}(\vk v_1,\vk v_2)=1}}\sum_{\substack{h_1\in {\tt G}_{\vk{v}_1}(H_i) \\ h_2\in {\tt G}_{\vk{v}_2}(H_j)}}\covv{\mathbb{I}_n(h_1,s),\,\mathbb{I}_n(h_2,t)}\notag\\
	&= \sum_{\vk v\in{\tt K}^{V_{i,j}^1}(\scale{n})}\sum_{\substack{\vk{v}_1\in{\tt K}^{\mathcal{V}(H_i)}(\vk{v})\\\vk{v}_2\in{\tt K}^{\mathcal{V}(H_j)}(\vk{v}) \\ {\tt I}(\vk v_1,\vk v_2)=1}}(\mathcal{V}(H_i)-1)!(\mathcal{V}(H_j)-1)!\mathcal{C}_{H_i,H_j}(s,t)\notag\\
	&= \sum_{\substack{\vk{v}_1\in{\tt K}^{\mathcal{V}(H_i)}(\scale{n})\\\vk{v}_2\in{\tt K}^{\mathcal{V}(H_j)}(\scale{n}) \\ {\tt I}(\vk v_1,\vk v_2)=1}}(\mathcal{V}(H_i)-1)!(\mathcal{V}(H_j)-1)!\mathcal{C}_{H_i,H_j}(s,t).\label{cov_repr_2}
\end{align}
Note that all terms of the sum \eqref{cov_repr_2} are equal, so we are left to count the number of terms involved in the sum. By straightforward calculation we obtain that
\begin{align}
	\sum_{\substack{\vk{v}_1\in{\tt K}^{\mathcal{V}(H_i)}(\scale{n})\\\vk{v}_2\in{\tt K}^{\mathcal{V}(H_j)}(\scale{n}) \\ {\tt I}(\vk v_1,\vk v_2)=1}} 1&=\left(n\frac{(n-1)!}{(\mathcal{V}(H_i)-1)!(n-\mathcal{V}(H_i))!}\frac{(n-\mathcal{V}(H_i))!}{(\mathcal{V}(H_j)-1)!(n-\mathcal{V}(H_i)-\mathcal{V}(H_j)+1)!}\right)\notag\\
	&= \frac{n!}{(n-\mathcal{V}(H_i)-\mathcal{V}(H_j)+1)!}\frac{1}{\bigl((\mathcal{V}(H_i)-1)!\bigr)\bigl((\mathcal{V}(H_j)-1)!\bigr)}.\label{cov_repr_3}
\end{align}
Combining \eqref{cov_repr_2} with \eqref{cov_repr_3} we obtain, as $n\to\infty$,
\begin{align}
	&\sum_{\vk h\in {\tt G}_{\scale{n}}^{i,j,1}}\covv{\mathbb{I}_n(h_1,s),\,\mathbb{I}_n(h_2,t)}\sim n^{\mathcal{V}(H_i)+\mathcal{V}(H_j)-1}\mathcal{C}_{H_i,H_j}(s,t).\label{cov_repr_4}
\end{align}
We are left to show that the last sum of \eqref{cov_repr_1} is asymptotically negligible comparing to the asymptotics obtained in \eqref{cov_repr_4}. We look at each term of this sum separately. The restriction $\abs{{\tt V}(h_1)\cap {\tt V}(h_2)}=k$ implies that $\abs{{\tt V}(h_1)\cup {\tt V}(h_2)} = V^k_{i,j}$. Using that the covariance of two Bernoulli random variables cannot be greater than 1, we get
\begin{align}
	\abs{\sum_{\vk h\in {\tt G}_{\scale{n}}^{i,j,k}}\covv{\mathbb{I}_n(h_1,s),\,\mathbb{I}_n(h_2,t)}}\leqslant \sum_{\vk h\in {\tt G}_{\scale{n}}^{i,j,k}} 1=\sum_{\vk v\in{\tt K}^{V^k_{i,j}}(\scale{n})}\sum_{\vk h\in {\tt G}_{\vk v}^{i,j,k}}1.\label{cov_repr_5}
\end{align}
Hence, it is enough to bound the number of terms in \eqref{cov_repr_5}. First, using \eqref{G_cardinality} we obtain (recall that $\vk v\in{\tt K}^{V^k_{i,j}}$, hence $\abs{\vk v} = V^k_{i,j}$)
\begin{align}
	\sum_{\vk h\in {\tt G}_{\vk v}^{i,j,k}}1&\leqslant \sum_{\vk h\in {\tt G}_{\vk v}^{i,j}}1 = \abs{{\tt G}_{\vk v}(H_i)}\cdot\abs{{\tt G}_{\vk v}(H_j)} = \frac{\bigl(V^k_{i,j}\bigr)!}{\left(V^k_{i,j}-\mathcal{V}(H_i)\right)!\mathcal{A}(H_i)}\frac{\bigl(V^k_{i,j}\bigr)!}{\left(V^k_{i,j}-\mathcal{V}(H_j)\right)!\mathcal{A}(H_j)}\notag\\
	&\leqslant\left(\bigl(V^k_{i,j}\bigr)!\right)^2.\label{cov_repr_6}
\end{align}
Then, by straightforward calculation
\begin{align} 
	\sum_{\vk v\in{\tt K}^{V^k_{i,j}}(\scale{n})}\left(\bigl(V^k_{i,j}\bigr)!\right)^2 = \frac{n!}{\bigl(n- V^k_{i,j}\bigr)!\bigl(V^k_{i,j}\bigr)!}\left(\bigl(V^k_{i,j}\bigr)!\right)^2\leqslant  \bigl(V^k_{i,j}\bigr)!n^{V^k_{i,j}}.\label{cov_repr_7}
\end{align}
Combining \eqref{cov_repr_5} with \eqref{cov_repr_6} and \eqref{cov_repr_7} we obtain, for any $k>1$ and as $n\to\infty$,
\begin{align}
	\abs{\sum_{\vk h\in {\tt G}_{\scale{n}}^{i,j,k}}\covv{\mathbb{I}_n(h_1,s),\,\mathbb{I}_n(h_2,t)}} \leqslant \bigl(V^k_{i,j}\bigr)!n^{V^k_{i,j}} =  o\bigl(n^{\mathcal{V}(H_i) + \mathcal{V}(H_j) - 1}\bigr).\label{cov_repr_8}
\end{align}
Hence, \eqref{cov_formula} follows after combining \eqref{cov_repr_1}, \eqref{cov_repr_4} and \eqref{cov_repr_8}. This concludes the proof of \nelem{exp_and_cov}.
\QED

\subsection{Proof of \nelem{asymptotic_normality}}\label{Section:proof_fdd_convergence}

To establish the claimed convergence, we use the Cram\'er--Wold device as follows. Fix some $0\leq t_1 <...<t_d\leq T$ for some $d\in\mathbb{N}$. It is then enough to show that for any vectors $\vk u_1,\ldots \vk u_d\in\R^m$, the random variable
\begin{align}
	\xi = \sum_{p=1}^{d}\scale{\vk u_p, \vk X^{\star}_n(t_p)} = \sum_{p=1}^{d}\sum_{i=1}^{m}u_{p,i}X^{\star}_{n,i}(t_p)\label{xi_init_def}
\end{align}
converges in distributions to a centered Gaussian random variable with the respective variance. 

Using \cite[Theorem 30.2]{billingsley2017probability}, it is sufficient to justify, for any  $z\in\N$, the moment convergence
\begin{align}
	\lim_{n\to\infty}\E{\xi^z} = 
	\begin{cases} 
		\sigma^z(z-1)!! \qquad &\text{if $z$ is even}, \\
		0 \qquad &\text{if $z$ is odd},
	\end{cases}\label{normality_claim}
\end{align}
with the constant 
\begin{align}
	\sigma^2 = \sum_{p,q=1}^{d}\sum_{i,j=1}^{m}u_{p,i}u_{q,j}\mathcal{C}_{H_i,H_j}(t_{p},t_{q}).\label{sigma_def}
\end{align}
From \eqref{X_def} we know that, for every $i\in\scale{m}$, the process $X_{n,i}^{\star}(\cdot)$ can be represented as a sum of centered indicators over all subgraphs $h\in{\tt G}_{\scale{n}}(H_i)$. We would like to do the same for $\xi$. For any $i\in\scale{m}$ and any $h\in{\tt G}_{\scale{n}}(H_i)$, define
\begin{align}
	\xi_i(h) = \sum_{p=1}^{d}u_{p,i}\bigl(\mathbb{I}_n(h,t_p) - \mathcal{P}_{H_i}(t_p)\bigr).\label{xi_def}
\end{align} 
Hence, according to \eqref{X_def} and \eqref{X_star_def} we can represent $\xi$ defined in \eqref{xi_init_def} as a sum over all possible $\xi_i(h)$ as
\begin{align}
	\xi = \sum_{i=1}^{m}\sum\limits_{h\in{\tt G}_{\scale{n}}(H_i)}\frac{\xi_i(h)}{n^{\mathcal{V}(H_i) - 1/2}}.\label{xi_sum_repr}
\end{align}

Taking expectation of the sum on the right hand side of \eqref{xi_sum_repr} to the power $z$ we obtain the following crucial expression of $z$'th moment of $\xi$: 
\begin{align}
	\E{\xi^z} &= \sum_{i_1,\ldots,i_z=1}^{m}\E{\prod_{k=1}^{z}\left(\sum_{h\in{\tt G}_{\scale{n}}(H_{i_k})}\frac{\xi_{i_k}(h)}{n^{\mathcal{V}(H_{i_k}) - 1/2}}\right)}\notag\\
	&=\sum_{\vk i\in\scale{m}^z}\sum_{\substack{h_1\in{\tt G}_{\scale{n}}(H_{i_1}) \\ \ldots \\ h_{z}\in{\tt G}_{\scale{n}}(H_{i_z})}}\E{\prod_{k=1}^{z}\frac{\xi_{i_k}(h_k)}{n^{\mathcal{V}(H_{i_k}) - 1/2}}}\notag\\
	&=\sum_{\vk i\in\scale{m}^z}\frac{1}{n^{\sum_{k=1}^{z}\mathcal{V}(H_{i_k}) - z/2}}\sum_{\vk h\in {\tt G}^{\vk i}_{\scale{n}}}\E{\prod_{k=1}^{z}\xi_{i_k}(h_k)},\label{moment_calculation_1}
\end{align}
where, to simplify the notations, for any ${\tt S}\subset\N$ and any $\vk i\in\scale{m}^z$ we denote
\begin{align*}
	{\tt G}^{\vk i}_{{\tt S}} = \bigtimes\limits_{k=1}^{z}{\tt G}_{{\tt S}}(H_{i_k}).
\end{align*}
We start the analysis of the right-hand side of \eqref{moment_calculation_1} by investigating, given $\vk i\in\scale{m}$, for which collections $\vk h\in{\tt G}^{\vk i}_{\scale{n}}$ the corresponding expectation is different from zero.

Using \eqref{I_N_def} we obtain that, if ${\tt V}(h)\cap{\tt V}(h^{\prime})=\varnothing$ for some $k,k^{\prime}\in\scale{z}$ and some subgraphs $h\in {\tt G}_{\scale{n}}(H_{i_k})$, $h^{\prime}\in{\tt G}_{\scale{n}}(H_{i_k{^{\prime}}})$, then the random variables $\xi_{i_k}(h)$ and $\xi_{i_{k^{\prime}}}(h^{\prime})$ are independent. Hence, if for some $k^{\prime}\in\scale{z}$ and all $k\in\scale{z}\setminus\{k^{\prime}\}$ we have that ${\tt V}(h_{k})\cap {\tt V}(h_{k^{\prime}}) = 0$, then
\begin{align*}
	\E{\prod_{k=1}^{z}\xi_{i_k}(h_k)} = \E{\prod_{\substack{k=1 \\ k\not= k^{\prime}}}^{z}\xi_{i_k}(h_k)}\E{\xi_{i_{k^{\prime}}}(h_{k^{\prime}})} = 0.
\end{align*}
Thus, we can consider only collections $\vk h\in{\tt G}^{\vk i}_{\scale{n}}$ which belong to the following set
\begin{align}
	\overline{{\tt G}}^{\vk i}_{\scale{n}} = \bigl\{\vk h\in {\tt G}^{\vk i}_{\scale{n}}\colon \forall k^{\prime}\in\scale{z} \,\exists k\in\scale{z}\setminus \{k^{\prime}\}\colon {\tt V}(h_k)\cap {\tt V}(h_{k^{\prime}})\geqslant 1\bigr\}.\label{g_condition}
\end{align}

Applying this observation, \eqref{moment_calculation_1} can be then expressed as
\begin{align}
	\E{\xi^z} = \sum_{\vk i\in\scale{m}^z}\frac{1}{n^{\sum_{k=1}^{z}\mathcal{V}(H_{i_k}) - z/2}}\sum_{\vk h\in \overline{{\tt G}}^{\vk i}_{\scale{n}}}\E{\prod_{k=1}^{z}\xi_{i_k}(h_k)}.\label{moment_calculation_2}
\end{align}

Consider first the inner sum of \eqref{moment_calculation_2}. For each fixed $\vk i\in\scale{m}^{z}$ we classify all the collections $\vk h\in\overline{{\tt G}}^{\vk{i}}_{\scale{n}}$ by the total number of vertices, i.e., by the cardinality of the set 
\begin{align}
	{\tt V}(\vk h) := \bigcup_{k=1}^{z}{\tt V}\left(h_k\right).\label{V_vk_def}
\end{align}
Thus, we can first choose $\abs{{\tt V}\left(\vk h\right)}$ out of $n$ vertices where the graphs $h_1,\ldots,h_z$ are located (this number is clearly $n$-dependent), consider a complete graph $K_{\abs{{\tt V}\left(\vk h\right)}}$ on these vertices, and after we choose $g_i$'s as subgraphs of the chosen complete graph $K_{\abs{{\tt V}\left(\vk h\right)}}$ (this number is $n$-independent). We denote by $V^{\max}_{\vk i}$ the maximal possible number $\abs{{\tt V}(\vk h)}$ for $\vk h\in\overline{{\tt G}}^{\vk i}_{\scale{n}}$, then

\begin{align}
	\sum_{\vk h\in\overline{{\tt G}}^{\vk i}_{\scale{n}}}\E{\prod_{k=1}^{z}\xi_{i_k}(h_k)} &= \sum_{c=0}^{V^{\max}_{\vk i}}\sum_{\substack{\vk h\in\overline{{\tt G}}^{\vk i}_{\scale{n}}\\ \abs{{\tt V}\left(\vk h\right)} = c }}\E{\prod_{k=1}^{z}\xi_{i_k}(h_k)}\notag\\
	&= \sum_{c=0}^{V_{\vk i}^{\max}}\sum_{\vk v\in{\tt K}^{c}(\scale{n})}\sum_{\substack{\vk h\in\overline{{\tt G}}^{\vk i}_{\vk v}\\ \abs{{\tt V}\left(\vk h\right)} = c }}\E{\prod_{k=1}^{z}\xi_{i_k}(h_k)}\notag\\
	&= \sum_{c=0}^{V_{\vk i}^{\max}}\frac{n!}{(n-c)!c!}\sum_{\substack{\vk h\in\overline{{\tt G}}^{\vk i}_{\scale{c}}\\ \abs{{\tt V}\left(\vk h\right)} = c }}\E{\prod_{k=1}^{z}\xi_{i_k}(h_k)}\notag\\
	&\sim \frac{n^{V_{i}^{\max}}}{V_{\vk i}^{\max}!}\sum_{\substack{\vk h\in\overline{{\tt G}}^{\vk i}_{\scale{V_{\vk i}^{\max}}} \\ \abs{{\tt V}(\vk h)} = V_{\vk i}^{\max} }}\E{\prod_{k=1}^{z}\xi_{i_k}(h_k)}\label{moment_calculation_3}
\end{align}
as $n\to\infty$. Note that the sim in \eqref{moment_calculation_3} does not depend on $n$, hence all the asymptotics is located in the factor $n^{V_{\vk i}^{\max}}$. So, the natural next step is to calculate $V_{\vk i}^{\max}$ for every $\vk i\in\scale{m}^z$. 

Define a simple labeled graph $\mathcal{G}(\vk h)$ on $z$ vertices as follows. An edge between vertices $k,k^{\prime}\in\scale{z}$ exists if and only if ${\tt V}(h_{k})\cup{\tt V}(h_{k^{\prime}})\not=\varnothing$, and let $\mathcal{G}^{\top}(\vk h)$ be a spanning forest of $\mathcal{G}(\vk h)$. Each edge $e_{k,k^{\prime}}$ of the forest $\mathcal{G}^{\top}(\vk h)$ means that the graphs $h_{k}$ and $h_{k^{\prime}}$ have at least one common vertex. Thus,
\begin{align*}
	\abs{{\tt V}\left(\vk h\right)}\leqslant \sum_{k=1}^{z}\mathcal{V}(H_{i_k}) - \mathcal{E}\bigl(\mathcal{G}^{\top}(\vk h)\bigr).
\end{align*}

Condition \eqref{g_condition} tells us that each vertex of the graph $\mathcal{G}(\vk h)$ has degree at least one, implying that
\begin{align*}
	\mathcal{E}\bigl(\mathcal{G}^{\top}(\vk h)\bigr) \geqslant \left\lceil\frac{z}{2}\right\rceil.
\end{align*}
In particular, when $z$ is odd and for any $\vk i\in\scale{m}^z$ such that $V_{\vk i}^{\max}< \sum_{k=1}^{z}\mathcal{V}(H_{i_k}) - z/2$, from \eqref{moment_calculation_3} we obtain, as $n\to\infty$,
\begin{align}
	\sum_{\vk h\in\overline{{\tt G}}^{\vk{i}}_{\scale{n}}}\E{\prod_{k=1}^{z}\xi_{i_k}(h_k)} = o\left(n^{\sum_{k=1}^{z}\mathcal{V}(H_{i_k}) - z/2}\right).\label{momant_calculation_z_odd}
\end{align}
Thus, after combining \eqref{moment_calculation_2} and \eqref{momant_calculation_z_odd}, \eqref{normality_claim} follows for odd $z$. Thus, we are left to consider only the case of even $z$.

When $z$ is even, we restrict our attention to those $\vk h\in{\tt G}^{\vk i}_{\scale{n}}$ for which $\mathcal{G}(\vk h)$ forms a perfect matching, i.e., every connected component of $\mathcal{G}(\vk h)$ consists of exactly two vertices (because $\mathcal{G}(\vk h)$ has no isolated vertices from \eqref{g_condition}). Note that $\mathcal{E}\bigl(\mathcal{G}^{\top}(\vk h)\bigr) = z/2$ if and only if $\mathcal{G}(\vk h)$ is a perfect matching. We denote by ${\tt M}_z$ all possible perfect matchings with $z$ vertices labeled $1,\ldots,z$. Note that in this case $ V_{\vk i}^{\max} = \sum_{k=1}^{z}\mathcal{V}(H_{i_k}) - z/2$. Under the condition \eqref{g_condition}, we obtain that $\abs{{\tt V}\left(\vk h\right)} = V_{\vk i}^{\max}$ if and only if $\mathcal{G}(\vk h)\in{\tt M}_{z}$. In addition, for every $k,k^{\prime}\in\scale{z}$, $k<k^{\prime}$, if $(k,k^{\prime})\in {\tt E}\bigl(\mathcal{G}(\vk h)\bigr)$ then $\abs{{\tt V}(h_k)\cap{\tt V}(h_{k^{\prime}})}=1$; otherwise, if $(k,k^{\prime})\not\in {\tt E}\bigl(\mathcal{G}(\vk h)\bigr)$, then ${\tt V}(h_k)\cap {\tt V}(h_{k^{\prime}}) = \varnothing$. Thus, for even $z$, from \eqref{moment_calculation_3} we obtain, as $n\to\infty$,
\begin{align}
	\sum_{\vk h\in\overline{{\tt G}}^{\vk i}_{\scale{n}}}\E{\prod_{k=1}^{z}\xi_{i_k}(h_k)} &\sim \frac{n^{V_{\vk i}^{\max}}}{V_{\vk i}^{\max}!}\sum_{\substack{\vk h\in\overline{{\tt G}}^{\vk i}_{\scale{V_{\vk i}^{\max}}} \\ \abs{{\tt V}(\vk h)} = V_{\vk i}^{\max} }}\E{\prod_{k=1}^{z}\xi_{i_k}(h_k)}\notag\\
	&= \frac{n^{V_{\vk i}^{\max}}}{V_{\vk i}^{\max}!}\sum_{M\in{\tt M}_z}\sum_{\substack{\vk h\in\overline{{\tt G}}^{\vk i}_{\scale{V_{\vk i}^{\max}}} \\ \abs{{\tt V}(\vk h)} = V_{\vk i}^{\max} \\ \mathcal{G}(\vk h) = M}}\E{\prod_{k=1}^{z}\xi_{i_k}(h_k)}\notag\\
	&= \frac{n^{V_{\vk i}^{\max}}}{V_{\vk i}^{\max}!}\sum_{M\in{\tt M}_z}\sum_{\substack{\vk h\in\overline{{\tt G}}^{\vk i}_{\scale{V_{\vk i}^{\max}}} \\ \abs{{\tt V}(\vk h)} = V_{\vk i}^{\max} \\ \mathcal{G}(\vk h) = M}}\prod_{(k,k^{\prime})\in{\tt E}(M)}\E{\xi_{i_k}(h_k)\xi_{i_{k^{\prime}}}(h_{k^{\prime}})}.\label{moment_calculation_4}
\end{align}

At this point, we proceed with a detailed analysis of the right-hand side of \eqref{moment_calculation_4}, firstly considering the inner product. We number the edges of the graph $M$ as ${\tt E}(M) = \{(a_1,b_1),\ldots,(a_{z/2},b_{z/2})\}$. Thus, by using \eqref{xi_def}, we get
\begin{align}
	\prod_{(k,k^{\prime})\in{\tt E}(M)}\E{\xi_{i_k}(h_k)\xi_{i_{k^{\prime}}}(h_{k^{\prime}})} &= \prod_{k=1}^{z/2}\E{\xi_{i_{a_k}}(h_{a_k})\xi_{i_{b_k}}(h_{b_k})}\notag\\
	&=\prod_{k=1}^{z/2}\mathbb{E}\left\{\left( \sum_{p=1}^{d}v_{p,i_{a_k}}\bigl(\mathbb{I}_n(h_{a_k},t_p) - \mathcal{P}_{H_{i_{a_k}}}(t_p)\bigr)\right)\right.\notag\\
	&\qquad\qquad\qquad\times\left.\left( \sum_{q=1}^{d}v_{q,i_{b_k}}\bigl(\mathbb{I}_n(h_{b_k},t_q) - \mathcal{P}_{H_{i_{b_k}}}(t_q)\bigr)\right)\right\}\notag\\
	&=\prod_{k=1}^{z/2}\left(\sum_{p,q=1}^{d}u_{p,i_{a_k}}u_{q,i_{b_k}}\operatorname{Cov}\bigl(\mathbb{I}_{n}(h_{a_k},t_p),\,\mathbb{I}_{n}(h_{b_k},t_q)\bigl)\right)\notag\\
	&= \sum_{\vk p,\vk q\in\scale{d}^{z/2}}\prod_{k=1}^{z/2}\left(u_{p_k,i_{a_k}}u_{q_k,i_{b_k}}\operatorname{Cov}\bigl(\mathbb{I}_{n}(h_{a_k},t_{p_k}),\,\mathbb{I}_{n}(h_{b_k},t_{q_k})\bigl)\right).\label{moment_calculation_5}
\end{align}
We simplify the right-hand side of \eqref{moment_calculation_5} by denoting for $n\in\N$, $M\in{\tt M}_z$, $\vk h\in\overline{{\tt G}}^{\vk i}_{\scale{V_{\vk i}^{\max}}}$ and $\vk p,\vk q\in\scale{d}^{z/2}$,
\begin{align}
	Q_n(M, \vk i, \vk h,\vk p,\vk q) = \prod_{k=1}^{z/2}\operatorname{Cov}\bigl(\mathbb{I}_{n}(h_{a_k},t_{p_k}),\,\mathbb{I}_{n}(h_{b_k},t_{q_k})\bigl).\label{Q_n_def}
\end{align}
Using this notation, and combining \eqref{moment_calculation_3} with \eqref{moment_calculation_4} and \eqref{moment_calculation_5}, we obtain that, as $n\to\infty$ and for even $z\in\N$,

\begin{align}
	\E{\xi^z}&\sim\sum_{\vk i\in\scale{m}^z}\sum_{M\in{\tt M}_z}\sum_{\substack{\vk h\in\overline{{\tt G}}^{\vk i}_{\scale{V_{\vk i}^{\max}}} \\ \abs{{\tt V}(\vk h)} = V_{\vk i}^{\max} \\ \mathcal{G}(\vk h) = M}}\sum_{\vk p,\vk q\in\scale{d}^{z/2}}\frac{\left(\prod_{k=1}^{z/2}u_{p_k,i_{a_k}}u_{q_k,i_{b_k}}\right)Q_n(M, \vk i,\vk h,\vk p,\vk q)}{V_{\vk i}^{\max}!}\notag\\
	&=\sum_{M\in{\tt M}_z}\sum_{\vk i\in\scale{m}^z}\sum_{\vk p,\vk q\in\scale{d}^{z/2}}\left(\prod_{k=1}^{z/2}u_{p_k,i_{a_k}}u_{q_k,i_{b_k}}\right)\sum_{\substack{\vk h\in\overline{{\tt G}}^{\vk i}_{\scale{V_{\vk i}^{\max}}} \\ \abs{{\tt V}(\vk h)} = V_{\vk i}^{\max} \\ \mathcal{G}(\vk h) = M}}\frac{Q_n(M, \vk i,\vk h,\vk p,\vk q)}{V_{\vk i}^{\max}!}.\label{moment_calculation_6}
\end{align}

Consider now the last sum of \eqref{moment_calculation_6} and fix some $\vk i\in\scale{m}^{z}$. The conditions $\mathcal{G}(\vk h)=M$ and $\abs{{\tt V}(\vk h)} = V_{\vk i}^{\max}$ tell us that, for any $k\in\scale{z/2}$, we have $\abs{{\tt V}(h_{a_k})\cup{\tt V}(h_{b_k})} = \mathcal{V}(H_{i_{a_k}}) + \mathcal{V}(H_{i_{b_k}}) - 1$ and, for different $k,k^{\prime}\in\scale{z/2}$, the sets ${\tt V}(h_{a_k})\cup{\tt V}(h_{b_k})$ and ${\tt V}(h_{a_{k^{\prime}}})\cup{\tt V}(h_{b_{k^{\prime}}})$ are disjoint. Hence, we can first separate the set of vertices $\scale{V_{\vk i}^{\max}}$ into disjoint sets $\vk v_1,\ldots,\vk v_{z/2}$ of cardinality $\abs{\vk v_k} = \mathcal{V}(H_{i_{a_k}}) + \mathcal{V}(H_{i_{b_k}}) - 1$, respectively. We can then choose $\vk h\in\overline{\tt G}^{\vk i}_{\scale{V_{\vk i}^{\max}}}$ in such a way that, for any $k\in\scale{z/2}$, we have ${\tt V}(h_{a_k}),{\tt V}(h_{b_k})\subset \vk v_k$ and $\abs{{\tt V}(h_{a_k})\cap{\tt V}(h_{b_k})}=1$, thereby yelding
\begin{align}
	\sum_{\substack{\vk h\in\overline{{\tt G}}^{\vk i}_{\scale{V_{\vk i}^{\max}}} \\ \abs{{\tt V}(\vk h)} = V_{\vk i}^{\max} \\ \mathcal{G}(\vk h) = M}}\frac{Q_n(M, \vk i,\vk h,\vk p,\vk q)}{V_{\vk i}^{\max}!} &= \sum_{\substack{\vk v_1,\ldots,\vk v_{z/2}\subset\scale{V_{\vk i}^{\max}} \\ \bigcup_{k=1}^{z/2}\vk v_k = \scale{V_{\vk i}^{\max}} \\ \abs{\vk v_k} = \mathcal{V}(H_{i_{a_k}}) + \mathcal{V}(H_{i_{b_k}}) - 1}}\frac{Q^{\star}(M,\vk i,\vk p,\vk q)}{V_{\vk i}^{\max}},\label{moment_calculation_7}
\end{align}
where (recall \eqref{Q_n_def} and \eqref{C_def})
\begin{align}
	Q^{\star}(M,\vk i,\vk p,\vk q)&:=\sum_{\substack{h_{a_1}\in {\tt G}_{\vk v_1}(H_{i_{a_1}}),\\ h_{b_1}\in {\tt G}_{\vk v_1}(H_{i_{b_1}}),\\\abs{{\tt V}(h_{a_1})\cap{\tt V}(h_{b_1})}=1}}\ldots \sum_{\substack{h_{a_{z/1}}\in {\tt G}_{\vk v_{z/2}}(H_{i_{a_{z/2}}}),\\ h_{b_{z/2}}\in {\tt G}_{\vk v_{z/2}}(H_{i_{b_{z/2}}}),\\\abs{{\tt V}(h_{a_{z/2}})\cap{\tt V}(h_{b_{z/2}})}=1}}Q_n(M, \vk i,\vk h,\vk p,\vk q)\notag\\
	&= \prod_{k=1}^{z/2}\left(\sum_{\substack{h_{a_k}\in {\tt G}_{\vk v_k}(H_{i_{a_k}}),\\ h_{b_k}\in {\tt G}_{\vk v_k}(H_{i_{b_k}}),\\\abs{{\tt V}(h_{a_k})\cap{\tt V}(h_{b_k})}=1}}\operatorname{Cov}\bigl(\mathbb{I}_{n}(h_{a_k},t_{p_k}),\,\mathbb{I}_{n}(h_{b_k},t_{q_k})\bigl)\right)\notag\\
	&= \prod_{k=1}^{z/2}\biggl(\bigl((\mathcal{V}(H_{i_{a_k}})-1)!\bigr)\bigl((\mathcal{V}(H_{i_{b_k}})-1)!\bigr)\mathcal{C}_{H_{i_{a_{k}}},H_{i_{b_k}}}(t_{p_k},t_{q_k})\biggr).\label{moment_calculation_8}
\end{align}
Combining now \eqref{moment_calculation_7} with \eqref{moment_calculation_8}, and using the fact that
\begin{align*}
	\sum_{\substack{\vk v_1,\ldots,\vk v_{z/2}\subset\scale{V_{\vk i}^{\max}} \\ \bigcup_{k=1}^{z/2}\vk v_k = \scale{V_{\vk i}^{\max}} \\ \abs{\vk v_k} = \mathcal{V}(H_{i_{a_k}}) + \mathcal{V}(H_{i_{b_k}}) - 1}}\,1 = \frac{V_{\vk i}^{\max}}{\prod\limits_{k=1}^{z/2}\bigl((\mathcal{V}(H_{i_{a_k}})-1)!\bigr)\bigl((\mathcal{V}(H_{i_{b_k}})-1)!\bigr)},
\end{align*}
we obtain that
\begin{align}
	\sum_{\substack{\vk h\in\overline{{\tt G}}^{\vk i}_{\scale{V_{\vk i}^{\max}}} \\ \abs{{\tt V}(\vk h)} = V_{\vk i}^{\max} \\ \mathcal{G}(\vk h) = M}}\frac{Q_n(M, \vk i,\vk h,\vk p,\vk q)}{V_{\vk i}^{\max}!} = \left(\prod_{k=1}^{z/2}\mathcal{C}_{H_{i_{a_{k}}},H_{i_{b_k}}}(t_{p_k},t_{q_k})\right).\label{moment_calculation_9}
\end{align}
After substituting \eqref{moment_calculation_9} into \eqref{moment_calculation_6}, we obtain, as $n\to\infty$ and for even $z\in\N$ (recall the definition of $\sigma^2$ given in \eqref{sigma_def}),
\begin{align*}
	\E{\xi^z}&\sim\sum_{M\in{\tt M}_z}\sum_{\vk i\in\scale{m}^z}\sum_{\vk p,\vk q\in\scale{d}^{z/2}}\left(\prod_{k=1}^{z/2}u_{p_k,i_{a_k}}u_{q_k,i_{b_k}}\right)\left(\prod_{k=1}^{z/2}\mathcal{C}_{H_{i_{a_{k}}},H_{i_{b_k}}}(t_{p_k},t_{q_k})\right)\\
	&=\sum_{M\in{\tt M}_z}\prod_{k=1}^{z/2}\left(\sum_{p,q=1}^{d}\sum_{i,j=1}^{m}u_{p,i}u_{q,j}\mathcal{C}_{H_i,H_j}(t_p,t_q)\right)=\sum_{M\in{\tt M}_z}\bigl(\sigma^{2}\bigr)^{z/2}=\abs{{\tt M}_z}\sigma^{z}.
\end{align*}

Thus, \eqref{normality_claim} follows for even $z$ from the fact that $\abs{{\tt M}_z} = (z-1)!!$. This completes the proof.
\QED

\subsection{Proof of \nelem{lem:tightness}}\label{Section:proof_of_tightness} Due to the definition of Euclidean norm, we can write
\begin{align*}
	\E{\abs{\vk X^{\star}_n(t)- \vk X^{\star}_n(s)}^2\abs{\vk X^{\star}_n(s) - \vk X^{\star}_n(r)}^2} = \sum_{i,j=1}^{n}\Delta_{i,j,n}(r,s,t),
\end{align*}
where for any $i,j\in\scale{m}$, $n\in\N$ and $0\leqslant r<s<t\leqslant T$,
\begin{align*}
	\Delta_{i,j,n}(r,s,t) = \E{\abs{X^{\star}_{n,i}(t)- X^{\star}_{n,i}(s)}^2\abs{X^{\star}_{n,j}(s)-X^{\star}_{n,j}(r)}^2}.
\end{align*}
Thus, to verify \eqref{tightness_claim} it is enough to show the for any $i,j\in\scale{m}$, $n\in\N$ and $0\leqslant r<s<t\leqslant T$,
\begin{align}
	\Delta_{i,j,n}(r,s,t)\leqslant F_{i,j}^2\left(\bigl(2V_{i,j}^{1}\bigr)!\right)^3\abs{t-r}^2,\label{tightness_claim_ij}
\end{align}
where $V_{i,j}^{1}$ is defined in \eqref{V_ij_def}.

For every $i\in \scale{m}$, $0\leqslant s<t\leqslant T$ and $h\in{\tt G}_{\scale{n}}(H_i)$, we denote the centered increment of the indicator function \eqref{I_N_def} by
\begin{align}
	\eta_{h,n,i}(s,t) &= \mathbb{I}_n(h,t) - \mathcal{P}_{H_i}(t) - \mathbb{I}_n(h,s) + \mathcal{P}_{H_i}(s)\label{nu_def}.
\end{align}
Combining \eqref{X_star_def} with \eqref{X_def}, we obtain that, for any $n\in\N$, $i\in\scale{m}$ and $0\leqslant s<t\leqslant T$, we can express the increments of the process $X^{\star}_{n,i}(\cdot)$ as 
\begin{align*}
	X^{\star}_{n,i}(t)-X^{\star}_{n,i}(s) = \frac{1}{n^{\mathcal{V}(H_i)-1/2}}\sum_{h\in{\tt G}_{\scale{n}}(H_i)}\eta_{h,n,i}(s,t).
\end{align*}
Hence, for any $n\in\N$, $i,j\in\scale{m}$ and $0\leqslant r<s<t\leqslant T$,
\begin{align}
	\Delta_{n,i,j}(r,s,t) &= \frac{1}{n^{2\mathcal{V}(H_i) + 2\mathcal{V}(H_j) -2}}\sum_{\substack{h_1,h_2\in {\tt G}_{\scale{n}}(H_i)\\ h_3,h_4\in {\tt G}_{\scale{n}}(H_j)}} \E{\eta_{h_1,n,i}(s,t)\eta_{h_2,n,i}(s,t)\eta_{h_3,n,j}(r,s)\eta_{h_4,n,j}(r,s)}\notag\\
	&= \frac{1}{n^{2V_{i,j}^{1}}}\sum_{\vk h\in{\tt G}_{\scale{n}}^{ii,jj}} \E{\prod_{k=1}^{4}\eta_{h_k,n,l_k}(x_k,y_k)},\label{tightness_sum}
\end{align} 
where, for any ${\tt S}\subset\N$,
\begin{align*}
	{\tt G}_{{\tt S}}^{ii,jj} &= {\tt G}_{{\tt S}}(H_i)\times {\tt G}_{{\tt S}}(H_i)\times {\tt G}_{{\tt S}}(H_j)\times {\tt G}_{{\tt S}}(H_j),
\end{align*}
and
\begin{align*}
	x_k &= \begin{cases} s,\qquad k\in\{1,2\},\\
		r,\qquad k\in\{3,4\},
	\end{cases}\qquad 
	y_k = \begin{cases} t,\qquad k\in\{1,2\},\\
		s,\qquad k\in\{3,4\}.
	\end{cases}\qquad l_k = \begin{cases} i,\qquad k\in\{1,2\},\\
		j,\qquad k\in\{3,4\}.
	\end{cases}
\end{align*}

The first step is to distinguish the zero terms of \eqref{tightness_sum} from the non-zero ones. Next, we derive an upper bound for each non-zero term in \eqref{tightness_sum}, and finally we show that the sum of the derived bounds does not exceed the right-hand side of \eqref{tightness_claim_ij}.

We identify the non-zero terms in the sum appearing in \eqref{tightness_sum}. From the definition \eqref{nu_def} we notice that the random variables $\eta_{h_k,n,l_k}(x_k,y_k)$ are centered, i.e., for any $n\in\N$, $k\in\scale{4}$, $\E{\eta_{h_k,n,l_k}(x_k,y_k)} = 0$, and for any $n\in\N$ and $k,k^{\prime}\in\scale{4}$, the random variables $\eta_{h_k,n,l_k}(x_k,y_k)$ and $\eta_{h_{k^{\prime}},n,l_{k^{\prime}}}(x_{k^{\prime}},y_{k^{\prime}})$ are independent if ${\tt V}(h_k)\cap {\tt V}(h_{k^{\prime}}) =\varnothing$.
Hence, if there exist $k^{\prime}\in\scale{4}$ such that for every $k\in\scale{4}$, $k\not= k^{\prime}$, we deduce that ${\tt V}(h_k)\cap {\tt V}(h_{k^{\prime}}) = \varnothing$, thereby yelding
\begin{align*}
	\E{\prod_{k=1}^{4}\eta_{h_k,n,l_k}(x_k,y_k)} = \E{\eta_{h_{k^{\prime}},n,l_{k^{\prime}}}(x_{k^{\prime}},y_{k^{\prime}})}\E{\prod_{\substack{k=1 \\ k\not= k^{\prime}}}^{4}\eta_{h_k,n,l_k}(x_k,y_k)} = 0.
\end{align*}
So, a term in the sum \eqref{tightness_sum} is non-zero only if it corresponds to some $\vk g\in{\tt G}_{\scale{n}}^{ii,jj}$ such that, for any $k\in\scale{4}$, there exists $k^{\prime}\in\scale{4}$, $k^{\prime}\not= k$, with ${\tt V}(h_k)\cap {\tt V}(h_{k^{\prime}})\not=\varnothing$. In particular, this restriction means that 
	$\abs{{\tt V}(\vk h)}\leqslant 2V_{i,j}^{1}$,
where ${\tt V}(\vk h)$ is defined in \eqref{V_vk_def} by taking $z=4$. Hence, neglecting the zero terms of \eqref{tightness_sum}, we can represent the function $\Delta_{n,i,j}(r,s,t)$ for all $n\in\N$, $i,j\in\scale{m}$ and $0\leqslant r<s<t\leqslant T$ as
\begin{align}
	\Delta_{n,i,j}(r,s,t) &= \frac{1}{n^{2V_{i,j}^{1}}}\sum_{\vk h\in{\tt G}_{\scale{n}}^{ii,jj}} \E{\prod_{k=1}^{4}\eta_{h_k,n,l_k}(x_k,y_k)} \notag\\
	&= \frac{1}{n^{2V_{i,j}^{1}}}\sum_{\substack{\vk h\in{\tt G}_{\scale{n}}^{ii,jj} \\ \abs{{\tt V}(\vk h)}\leqslant 2V_{i,j}^{1}}} \E{\prod_{k=1}^{4}\eta_{h_k,n,l_k}(x_k,y_k)}.\label{tightness_sum_non_zero}
\end{align}

Next, we proceed to the second step, where we present an upper bound for each term of \eqref{tightness_sum_non_zero} separately. To this end, fix some $n\in\N$, $i,j\in\scale{m}$ and $\vk h\in{\tt G}_{\scale{n}}^{ii,jj}$. Define for any vertex $v\in\scale{n}$ and any $a,b\in[0,T]$ the event $\Omega_{v}(a,b)$, which tells us that vertex $v$ changed its opinion at least once between time points $a$ and $b$. Additionally, for any edge $(u,v)\in\scale{n}^2$ and any $a,b\in[0,T]$, we define the event $\Omega_{(u,v)}(a,b)$ which tells us that the Poisson clock attached to the edge between vertices $u$ and $v$ rings at least once. According to the setup of the model, all the events $\Omega_{v}(a,b)$, $\Omega_{(u,v)}(a,b)$ are independent of each other for different lower subscripts. For any $u,v\in\scale{n}$ and any $a,b\in[0,T]$, we can derive the following bounds for their probabilities:
\begin{align}
	\pk{\Omega_{v}(a,b)}\leqslant C\abs{a-b},\label{vertex_bound}\\
	\pk{\Omega_{(u,v)}(a,b)}\leqslant C\abs{a-b}\label{edge_bound},
\end{align}
where
	$C=2\max(\gamma_{-+},\gamma_{+-},\pi_+,\pi_-)$.
Thus, for any $k\in\scale{4}$, we can write
\begin{align*}
	\pk{\mathbb{I}_n(h_k,y_k)\not=\mathbb{I}_n(h_k,x_k)}&\leqslant \sum_{v\in{\tt V}(h_k)}\pk{\Omega_{v}(x_k,y_k)}+ \sum_{(u,v)\in{\tt E}(h_k)}\pk{\Omega_{(u,v)}(x_k,y_k)}\\
	&\leqslant C(\mathcal{V}(H_{l_k})+\mathcal{E}(H_{l_k}))\abs{y_k-x_k}.
\end{align*}
Additionally, for any $k\in\scale{4}$, on the event $\mathbb{I}_{n}(h_k,x_k) = \mathbb{I}_{n}(h_k,y_k)$ we obtain the following bound which holds everywhere on the aforementioned event
\begin{align}
	\abs{\eta_{h_k,n,l_k}(x_k,y_k)} &= \abs{\mathcal{P}_{H_{l_k}}(y_k) - \mathcal{P}_{H_{l_k}}(x_k)}\notag\\
	&= \lvert\pk{\mathbb{I}_{n}(h_k,y_k) = 1, \mathbb{I}_{n}(h_k,x_k)=0}- \pk{\mathbb{I}_{n}(h_k,y_k) = 0, \mathbb{I}_{n}(h_k,x_k)=1}\rvert\notag\\
	&\leqslant \pk{\mathbb{I}_{n}(h_k,y_k)\not=\mathbb{I}_{n}(h_k,x_k)}\leqslant C(\mathcal{V}(H_{l_k})+\mathcal{E}(H_{l_k}))\abs{t-r}.\label{eta_bound}
\end{align}

However, if $\mathbb{I}_{n}\{h_k,x_k) \not= \mathbb{I}_{n}(h_k,y_k)$,
then there is no appropriate bound for $\abs{\eta_{h_k,n,l_k}(x_k,y_k)}$, so we should find an upper bound for the probability of such event. To this end, we define two events, one of them dedicated to $h_1$ and $h_2$, and the other one to $h_3$ and $h_4$, as follows:
\begin{align*}
	\Omega_1 &= \bigl\{\mathbb{I}_{n}(h_1,y_1) \not= \mathbb{I}_{n}(h_1,x_1)\bigr\}\cup \bigl\{\mathbb{I}_{n}(h_2,y_2) \not= \mathbb{I}_{n}(h_2,x_2)\bigr\},\\
	\Omega_2 &= \bigl\{\mathbb{I}_{n}(h_3,y_3) \not= \mathbb{I}_{n}(h_3,x_3)\bigr\}\cup \bigl\{\mathbb{I}_{n}(h_4,y_4) \not= \mathbb{I}_{n}(h_4,x_4)\bigr\}.
\end{align*}
Then on the event $\Omega_1^{\rm c}$, which tells us that $\mathbb{I}_{n}(h_1,y_1) = \mathbb{I}_{n}(h_1,x_1)$ and $\mathbb{I}_{n}(h_2,y_2) = \mathbb{I}_{n}(h_2,x_2)$, we can apply \eqref{eta_bound} for both $\abs{\eta_{h_1,n,i}(x_1,y_1)}$ and $\abs{\eta_{h_2,n,i}(x_2,y_2)}$ obtaining the following inequality which holds everywhere on $\Omega_1^{\rm c}$
\begin{align}
	\abs{\eta_{h_1,n,i}(x_1,y_1)\eta_{h_2,n,i}(x_2,y_2)} \leqslant C^2 \bigl(\mathcal{V}(H_i) + \mathcal{E}(H_i)\bigr)^2\abs{t-r}^2,\label{O1_c}
\end{align}
similarly, on the event $\Omega_2^{\rm c}$ holds the following inequality
\begin{align}
	\abs{\eta_{h_3,n,j}(x_3,y_3)\eta_{h_4,n,j}(x_4,y_4)} \leqslant C^2 \bigl(\mathcal{V}(H_j) + \mathcal{E}(H_j)\bigr)^2\abs{t-r}^2.\label{O2_c}
\end{align}
Our goal is now to find upper bounds for the probabilities $\pk{\Omega_1}$, $\pk{\Omega_2}$ and $\pk{\Omega_1\Omega_2}$. Note that $\Omega_1$ implies that either at least one vertex of $h_1$ or $h_2$, or at least one edge of $h_1$ or $h_2$ changed its status between time points $s$ and $t$, i.e.,
\begin{align}
	\Omega_1 \subset\bigcup_{v\in{\tt V}(h_1)\cup {\tt V}(h_2)}\Omega_{v}(s,t) \cup \bigcup_{(u,v)\in{\tt E}(h_1)\cup {\tt E}(h_2)}\Omega_{(u,v)}(s,t).\label{Omega_incl}
\end{align}
Hence, applying \eqref{vertex_bound} and \eqref{edge_bound}, we derive the following bound for $\pk{\Omega_1}$ (recall that $r<s<t$)
\begin{align*}
	\pk{\Omega_1} &\leqslant \sum_{v\in{\tt V}(h_1)\cup {\tt V}(h_2)}\pk{\Omega_{v}(s,t)}+\sum_{(u,v)\in{\tt E}(h_1)\cup {\tt E}(h_2)}\pk{\Omega_{(u,v)}(s,t)}\notag\\
	&\leqslant 2C\bigl(\mathcal{V}(H_i) + \mathcal{E}(H_i)\bigr)\abs{t-s}\leqslant 2C\bigl(\mathcal{V}(H_i) + \mathcal{E}(H_i)\bigr)\abs{t-r}.
\end{align*}
By similar arguments, we get
\begin{align*}
	\pk{\Omega_2}\leqslant 2C\bigl(\mathcal{V}(H_j) + \mathcal{E}(H_j)\bigr)\abs{t-r}.
\end{align*}
For the intersection $\Omega_1\Omega_2$, by simple calculations and using \eqref{Omega_incl}, we can see that
\begin{align*}
	\Omega_1\Omega_2&\subset \left(\bigcup_{v\in{\tt V}(h_1)\cup {\tt V}(h_2)}\Omega_{v}(s,t) \cup \bigcup_{(u,v)\in{\tt E}(h_1)\cup {\tt E}(h_2)}\Omega_{(u,v)}(s,t)\right)\\
	&\qquad \bigcap \left(\bigcup_{v\in{\tt V}(h_3)\cup {\tt V}(h_4)}\Omega_{v}(r,s) \cup \bigcup_{(u,v)\in{\tt E}(h_3)\cup {\tt E}(h_4)}\Omega_{(u,v)}(r,s)\right)\\
	&=\left(\bigcup_{v\in{\tt V}(h_1)\cup {\tt V}(h_2)}\Omega_{v}(s,t) \cap \bigcup_{v\in{\tt V}(h_3)\cup {\tt V}(h_4)}\Omega_{v}(r,s)\right)\\
	&\qquad\bigcup \left(\bigcup_{v\in{\tt V}(h_1)\cup {\tt V}(h_2)}\Omega_{v}(s,t)\cap \bigcup_{(u,v)\in{\tt E}(h_3)\cup {\tt E}(h_4)}\Omega_{(u,v)}(r,s)\right)\\
	&\qquad\bigcup \left(\bigcup_{(u,v)\in{\tt E}(h_1)\cup {\tt E}(h_2)}\Omega_{(u,v)}(s,t) \cap \bigcup_{v\in{\tt V}(h_3)\cup {\tt V}(h_4)}\Omega_{v}(r,s)\right)\\
	&\qquad\bigcup \left(\bigcup_{(u,v)\in{\tt E}(h_1)\cup {\tt E}(h_2)}\Omega_{(u,v)}(s,t) \cap \bigcup_{(u,v)\in{\tt E}(h_3)\cup {\tt E}(h_4)}\Omega_{(u,v)}(r,s)\right)\\
	&=\left(\bigcup_{\substack{v_1\in{\tt V}(h_1)\cup {\tt V}(h_2) \\  v_2\in{\tt V}(h_3)\cup {\tt V}(h_4)}}\Omega_{v_1}(s,t)  \Omega_{v_2}(r,s)\right) \bigcup  \left(\bigcup_{\substack{v_1\in{\tt V}(h_1)\cup {\tt V}(h_2) \\ (u_2,v_2)\in{\tt E}(h_3)\cup {\tt E}(h_4)}}\Omega_{v_1}(s,t) \Omega_{(u_2,v_2)}(r,s)\right)\\
	&\qquad \bigcup \left(\bigcup_{\substack{(u_1,v_1)\in{\tt E}(h_1)\cup {\tt E}(h_2) \\ v_2\in{\tt V}(h_3)\cup {\tt V}(h_4)}}\Omega_{(u_1,v_1)}(s,t)  \Omega_{v_2}(r,s)\right)\\
	&\qquad \bigcup \left(\bigcup_{\substack{(u_1,v_1)\in{\tt E}(h_1)\cup {\tt E}(h_2) \\ (u_2,v_2)\in{\tt E}(h_3)\cup {\tt E}(h_4)}}\Omega_{(u_1,v_1)}(s,t)  \Omega_{(u_2,v_2)}(r,s)\right)
\end{align*}
Hence, applying \eqref{vertex_bound} and \eqref{edge_bound}, for $\pk{\Omega_1\Omega_2}$ we can derive the following upper bound
\begin{align}
	\pk{\Omega_1\Omega_2}&\leqslant \sum_{\substack{v_1\in{\tt V}(h_1)\cup {\tt V}(h_2) \\  v_2\in{\tt V}(h_3)\cup {\tt V}(h_4)}}\pk{\Omega_{v_1}(s,t)  \Omega_{v_2}(r,s)} +\sum_{\substack{v_1\in{\tt V}(h_1)\cup {\tt V}(h_2) \\ (u_2,v_2)\in{\tt E}(h_3)\cup {\tt E}(h_4)}}\pk{\Omega_{v_1}(s,t) \Omega_{(u_2,v_2)}(r,s)}\notag\\
	&\qquad + \sum_{\substack{(u_1,v_1)\in{\tt E}(h_1)\cup {\tt E}(h_2) \\ v_2\in{\tt V}(h_3)\cup {\tt V}(h_4)}}\pk{\Omega_{(u_1,v_1)}(s,t)  \Omega_{v_2}(r,s)}\notag\\
	&\qquad + \sum_{\substack{(u_1,v_1)\in{\tt E}(h_1)\cup {\tt E}(h_2) \\ (u_2,v_2)\in{\tt E}(h_3)\cup {\tt E}(h_4)}}\pk{\Omega_{(u_1,v_1)}(s,t)  \Omega_{(u_2,v_2)}(r,s)}\notag\\
	&= \sum_{\substack{v_1\in{\tt V}(h_1)\cup {\tt V}(h_2) \\  v_2\in{\tt V}(h_3)\cup {\tt V}(h_4)}}\pk{\Omega_{v_1}(s,t)}\pk{\Omega_{v_2}(r,s)} \notag\\
	&\qquad +\sum_{\substack{v_1\in{\tt V}(h_1)\cup {\tt V}(h_2) \\ (u_2,v_2)\in{\tt E}(h_3)\cup {\tt E}(h_4)}}\pk{\Omega_{v_1}(s,t)}\pk{\Omega_{(u_2,v_2)}(r,s)}\notag\\
	&\qquad + \sum_{\substack{(u_1,v_1)\in{\tt E}(h_1)\cup {\tt E}(h_2) \\ v_2\in{\tt V}(h_3)\cup {\tt V}(h_4)}}\pk{\Omega_{(u_1,v_1)}(s,t)}\pk{\Omega_{v_2}(r,s)}\notag\\
	&\qquad + \sum_{\substack{(u_1,v_1)\in{\tt E}(h_1)\cup {\tt E}(h_2) \\ (u_2,v_2)\in{\tt E}(h_3)\cup {\tt E}(h_4)}}\pk{\Omega_{(u_1,v_1)}(s,t)}\pk{\Omega_{(u_2,v_2)}(r,s)}\notag\\
	&\leqslant 4C^2\left(\mathcal{V}(H_i) + \mathcal{E}(G_i)\right)\left(\mathcal{V}(H_j) + \mathcal{E}(G_j)\right)\abs{t-r}^2.\label{O12}
\end{align}

Combining \eqref{O1_c}, \eqref{O2_c} and \eqref{O12} together with the bound $\eta_{h_k,n,l_k}(x_k,y_k)\leqslant 2$, we get
\begin{align}
	\E{\prod_{k=1}^{4}\eta_{h_k,n,l_k}(x_k,y_k)} &= \E{\prod_{i=1}^{4}\eta_{h_k,n,l_k}(x_i,y_i)\mathbb{I}\{\Omega_1\}\mathbb{I}\{\Omega_2\}} \notag\\
	&\quad + \E{\prod_{k=1}^{4}\eta_{h_k,n,l_k}(x_k,y_k)\mathbb{I}\{\Omega_1\}\mathbb{I}\{\Omega_2^{\rm c}\}} \notag\\ &\quad+\E{\prod_{k=1}^{4}\eta_{h_k,n,l_k}(x_k,y_k)\mathbb{I}\{\Omega_1^{\rm c}\}}\notag \\
	&\leqslant 16\pk{\Omega_1\Omega_2} + 4C^2\bigl(\mathcal{V}(H_j)+\mathcal{E}(H_j)\bigr)^2\pk{\Omega_1}\notag\\
	&\qquad + 4C^2\bigl(\mathcal{V}(H_i)+\mathcal{E}(H_i)\bigr)^2\pk{\Omega_1}\notag\\
	&à\leqslant F^2_{i,j}\abs{t-r}^2, \label{prod_eta_bound}
\end{align}
where the constant $F_{i,j}$ is defined in \eqref{F_ij_def}.

The inequality in \eqref{prod_eta_bound} provides upper bounds for each term of \eqref{tightness_sum_non_zero}. Hence, applying this bound in \eqref{prod_eta_bound} for each term of \eqref{tightness_sum_non_zero} and using \eqref{G_cardinality}, for any $n\in\N$, $i,j\in\scale{m}$ and $0\leqslant r<s<t\leqslant T$, we obtain 
\begin{align*}
	\Delta_{n,i,j}(r,s,t) &= \frac{1}{n^{2V_{i,j}^{1}}} \sum_{\substack{\vk h\in{\tt G}_{\scale{n}}^{ii,jj} \\ \abs{{\tt V}(\vk h)}\leqslant 2V_{i,j}^{1}}} \E{\prod_{k=1}^{4}\eta_{h_k,n,l_k}(x_k,y_k)}\\
	&\leqslant \frac{F^2_{i,j}(t-r)^2}{n^{2V_{i,j}^{1}}} \sum_{\vk v\in{\tt K}^{2V_{i,j}^{1}}(\scale{n})}\sum_{\vk h\in {\tt G}^{ii,jj}_{\vk v}} 1\\
	&= \frac{F^2_{i,j}(t-r)^2}{n^{2V_{i,j}^{1}}}\sum_{\vk v\in{\tt K}^{2V_{i,j}^{1}}(\scale{n})} \left(\frac{\bigl(2V_{i,j}^{1}\bigr)!}{(2V_{i,j}^{1}-\mathcal{V}(H_i))!\mathcal{A}(H_i)}\right)^2\left(\frac{\bigl(2V_{i,j}^{1}\bigr)!}{(2V_{i,j}^{1}-\mathcal{V}(H_j))!\mathcal{A}(H_j)}\right)^2\\
	&\leqslant F_{i,j}^2(t-r)^2 \left(\bigl(2V_{i,j}^{1}\bigr)!\right)^3.
\end{align*}

Thus, \eqref{tightness_claim_ij} follows. This concludes the proof.
\QED

\subsection{Proof of \neprop{prop:variance}}\label{Section:proofprop}

Given any $\vk v_1\in{\tt K}^{\mathcal{V}(H)}(\scale{n})$, $\vk v_2\in{\tt K}^{\mathcal{V}(H^{\prime})}(\scale{n})$ such that $\vk v_1\cap \vk v_2 = \{k\}$, fix $h\in {\tt G}_{\vk v_1}(H)$ and $h^{\prime}\in {\tt G}_{\vk v_2}(H^{\prime})$. Then, for any $t\in[0,T]$,
\[
\begin{array}{ll}
	\displaystyle \E{\I_n(h,t)\I_n(h^{\prime},t)} \\
	\qquad = \displaystyle\E{\prod_{\substack{(u,v)\in{\tt E}(h) \\ (w,z)\in{\tt E}(h^{\prime})}}a_{u,v}(t)a_{w,z}(t)\prod_{\substack{i\in{\tt V}(h) \\ j\in{\tt V}(h^{\prime})}}\mathbb{I}\{x_{i}(G_n(t)) = x_{i}(h)\}\mathbb{I}\{x_{j}(G_n(t)) = x_{j}(h^{\prime})\}}.
\end{array}
\]
Since the events involving the presence of active edges and the opinions of the vertices at time $t$ when considering both $h$ and $h^{\prime}$ are independent after conditioning on the type of the common vertex $k$, we deduce that we can write the above formula as 
\begin{equation}
	\int_0^1 \mathbb{E}_1 \ \mathbb{E}_2 \ \mathbb{P}(y_k(t)=\dd \ell_k),
\end{equation}
where
\begin{equation}\label{eq:term1}
	\begin{array}{ll}
		\mathbb{E}_1 &= \displaystyle\E{\prod_{\substack{(u,v)\in{\tt E}(h) \\ (w,z)\in{\tt E}(h^{\prime})}}a_{u,v}(t)a_{w,z}(t) | y_k(t) = \dd \ell_k} \\
		&= \displaystyle\E{\prod_{(u,v)\in{\tt E}(h)} a_{u,v}(t) | y_k(t) = \dd \ell_k} \E{\prod_{ (w,z)\in{\tt E}(h^{\prime})}a_{w,z}(t) | y_k(t) = \dd \ell_k}
	\end{array}
\end{equation}
and
\begin{equation}\label{eq:term2}
	\mathbb{E}_2 = \displaystyle\E{\prod_{\substack{i\in{\tt V}(h) \\ j\in{\tt V}(h^{\prime})}}\mathbb{I}\{x_{i}(G_n(t)) = x_{i}(h)\}\mathbb{I}\{x_{j}(G_n(t))= x_{j}(h^{\prime})\}| y_k(t) = \dd \ell_k}.
\end{equation}
Since the events involving the presence of active edges are mutually independent conditionally on the type of each vertex, and the evolution of the opinions of each vertex is independent from the others, note that \eqref{eq:term1} becomes
\begin{equation}\label{eq:term3}
	\begin{array}{ll}
		\displaystyle \int_{[0,1]^{|h|-1}} \prod_{(u,v)\in{\tt E}(h)} \mathcal{H}(t;\ell_u,\ell_v) \prod_{i\in {\tt V}(h)\setminus\{k\}} \mathbb{P}(y_i = \dd \ell_i) \\
		\qquad \qquad \qquad \qquad \qquad \times \displaystyle\int_{[0,1]^{|h^{\prime}|-1}} \prod_{(w,z)\in{\tt E}(h^{\prime})} \mathcal{H}(t;\ell_w,\ell_z) \prod_{i\in {\tt V}(h^{\prime})\setminus\{k\}} \mathbb{P}(y_i = \dd \ell_i).
	\end{array}
\end{equation}
Similarly, by nothing that the events involving the opinion of the vertices are mutually independent conditionally on the opinion of the common vertex $k$, and that the evolution of the opinions of each vertex is independent from the others, equation \eqref{eq:term2} becomes
\begin{equation}
	\begin{array}{ll}
		\displaystyle \sum_{o\in\{-,+\}} \E{\prod_{i\in{\tt V}(h)} \mathbb{I}\{x_{i}(G_n(t)) = x_{i}(h)\}| x_k(G_n(t))=o, y_k(t) = \dd \ell_k} \\
		\qquad \qquad \quad \displaystyle \times \E{\prod_{j\in{\tt V}(h^{\prime})}\mathbb{I}\{x_{j}(G_n(t))= x_{j}(h^{\prime})\}| x_k(G_n(t))=o, y_k(t) = \dd \ell_k} \mathbb{P}(x_k(G_n(t))=o) \\
		= \displaystyle \sum_{o\in\{-,+\}} \prod_{i\in{\tt V}(h)} \mathbb{P}(\mathbb{I}\{x_{i}(G_n(t)) = x_{i}(h)\}| x_k(G_n(t))=o) \\
		\qquad \qquad \qquad \qquad \qquad \displaystyle \times \prod_{j\in{\tt V}(h^{\prime})} \mathbb{P}(\mathbb{I}\{x_{j}(G_n(t))= x_{j}(h^{\prime})\}| x_k(G_n(t))=o) \mathbb{P}(x_k(G_n(t))=o).
	\end{array}
\end{equation}

By using the fact that $\mathbb{P}(x_i(G_n(t))=+, y_i(t)=\dd \ell_i)=f_+(t,\ell_i)\dd \ell_i$ and $\mathbb{P}(x_i(G_n(t))=-, y_i(t)=\dd \ell_i)=f_-(t,\ell_i)\dd \ell_i$ for any vertex $i$, and after arguing similarly for the terms $\E{\I_n(h,t)}$ and $\E{\I_n(h^{\prime},t)}$, \eqref{eq:C_G(t,t)} follows.

Finally, we show that $C_H(t,t)>0$ for any $t\in[0,T]$ as follows. Let $\vk v_1\cap\vk v_2 = \{k\}$. For any $h\in {\tt G}_{\vk v_1}(H) \cup {\tt G}_{\vk v_2}(H)$, we define  a $(x_k(G_n(s)))_{s\in[0,t]}$-measurable random variable as
\begin{align*}
	\xi(h,t) = \E{\mathbb{I}_n(h,t)\mid \bigl(x_k(G_n(s))\bigr)_{s\in[0,t]}}.
\end{align*}
Then, for any $h\in {\tt G}_{\vk v_1}(H)$ and $h^{\prime}\in {\tt G}_{\vk v_1}(H)$, (recall \eqref{P_def})
\begin{align*}
	\operatorname{Cov}\bigl(\mathbb{I}_{n}(h,t),\,\mathbb{I}_{n}(h^{\prime},t)\bigl) &= \E{\mathbb{I}_{n}(h,t)\mathbb{I}_{n}(h^{\prime},t)} - \mathcal{P}_H^2(t)\\
	&=\mathbb{E}\biggl\{\E{\mathbb{I}_{n}(h,t)\mathbb{I}_{n}(h^{\prime},t)\mid \bigl(x_v(G_n(s))\bigr)_{s\in[0,t]}}\biggr\} - \mathcal{P}_H^2(t)\\
	&=\E{\xi(h,t)\xi(h^{\prime},t)} - \mathcal{P}_H^2(t).
\end{align*}
Hence,
\begin{align*}
	\sum_{\substack{h\in {\tt G}_{\vk v_1}(H) \\ h^{\prime}\in {\tt G}_{\vk v_2}(H)}}\operatorname{Cov}\bigl(\mathbb{I}_{n}(h,t),\,\mathbb{I}_{n}(h^{\prime},t)\bigl) &= \E{\sum_{\substack{h\in {\tt G}_{\vk v_1}(H) \\ h^{\prime}\in {\tt G}_{\vk v_2}(H)}}\xi(h,t)\xi(h^{\prime},t)} - \sum_{\substack{h\in {\tt G}_{\vk v_1}(H) \\ h^{\prime}\in {\tt G}_{\vk v_2}(H)}}\mathcal{P}_H^2(t)\\
	&=\operatorname{Cov}(\xi(t),\xi^{\prime}(t)),
\end{align*}
where
\begin{align*}
	\xi(t) = \sum_{h\in {\tt G}_{\vk v_1}(H)}\xi(h,t), \qquad \xi^{\prime}(t) = \sum_{h^{\prime}\in {\tt G}_{\vk v_2}(H)}\xi(h^{\prime},t)
\end{align*}
are two identically distributed $(x_k(G_n(s)))_{s\in[0,t]}$-measurable random variables. From the symmetry of our model, the $\sigma$-algebras generated by these two random variables are equal, and therefore $\xi(t) = \xi^{\prime}(t)$ almost surely for any $t\in[0,T]$. Hence, for any $t\in[0,T]$,
\begin{align*} 
	\operatorname{Cov}(\xi(t),\xi^{\prime}(t)) = \operatorname{Var}(\xi(t))>0,
\end{align*}
as $\xi(t)$ is a non-degenerate random variable. Hence, the claim follows from \eqref{C_def}.
\QED

\section{Two-way feedback model: some insights}\label{Section:two-ways_feedback}

In this section, we extend the one-way feedback model described in Section \ref{sec:defmodel} to a co-evolutionary model. The dynamics of the edges is the same as for the previous model, but the dynamics of the vertices is different. The model can be described as follows, where for completeness we repeat the description of the dynamics of the edges:

\smallskip
{\it Vertex dynamics.} Each vertex holds opinion $+$ or $-$, and is assigned an independent rate-$\beta$ Poisson clock. Each time the clock rings the vertex selects one of its neighbors uniformly at random and copies the opinion of that vertex.

\smallskip
{\it Edges dynamics.} Each edge is re-sampled at rate $1$, i.e., a rate-$1$ Poisson clock is attached to each edge and when the clock rings the edge is active with a probability that depends on the current opinion of the two connected vertices: with probability $\pi_+$ if the two vertices hold opinion $+$, with probability $\pi_-$ if the two vertices hold opinion $-$, and otherwise with probability $\tfrac12(\pi_++\pi_-)$.

\smallskip
Note that this network is indeed co-evolutionary: the opinions of the adjacent vertices affect the probability that an edge is active (through $\pi_+ $ and $\pi_-$), and the state of adjacent edges affect the opinions of the vertices (because the vertices copy the opinions of their neighbors).

\smallskip
However, in this model the calculation of the limiting distribution of the subgraph count process $X_n(\cdot)$ described in Section \ref{Section:main_results}, even marginally, is much more complex as explained in what follows. Let us fix some $t\in[0,T]$ and some deterministic voter graph $H$. For any two \textit{disjoint} vectors $\vk v_1,\vk v_2\in{\tt K}^{\mathcal{V}(H)}(\scale{n})$, we define the counterpart of the constant \eqref{C_def} for this two-way feedback model  as
\begin{align*}
	\mathcal{C}^{\prime}_{n,H}(t) &= \frac{1}{\bigl(\mathcal{V}(H)!\bigr)^2}\sum_{\substack{h\in {\tt G}_{\vk v_1}(H) \\ h^{\prime}\in {\tt G}_{\vk v_2}(H)}}\operatorname{Cov}\bigl(\mathbb{I}_{n}(h,t),\,\mathbb{I}_{n}(h^{\prime},t)\bigl) = \frac{\operatorname{Cov}\bigl(\mathbb{I}_{n}(h,t),\,\mathbb{I}_{n}(h^{\prime},t)\bigl)}{\bigl(\mathcal{A}(H)\bigr)^2},
\end{align*}
where we define ${\tt G}_{\vk v}(H)$ and $\mathbb{I}_{n}(h,t)$ as in \eqref{tt_G_def} and \eqref{I_N_def}, respectively, for any choice $h\in {\tt G}_{\vk v_1}(H)$ and $h^{\prime}\in {\tt G}_{\vk v_2}(H)$ due to symmetry of the model and \eqref{G_cardinality}. Notice that in the one-way feedback model such a constant simply equals zero, because in that model, if graphs $h$ and $h^{\prime}$ does not intersect, the respective random variables $\mathbb{I}_{n}(h,t)$ and $\mathbb{I}_{n}(h^{\prime},t)$ are independent. However, this is not true for the two-way feedback model due to the intricate relation between the joint dynamics of vertices and edges. One can observe that $\lim_{n\to\infty}\mathcal{C}^{\prime}_{H}(t)\in(0,\infty)$ for the two-way feedback model (see Table \ref{tab:simulation} for some numerical results when $n=100$ and $n=200$). Then, by applying the same arguments given in \nelem{exp_and_cov}, we can figure out that, as $n\to\infty$,
\begin{align*}
	\E{X_{n}(t)} &= \frac{n! \mathcal{P}_{H}(t)}{(n-\mathcal{V}(H))!\mathcal{A}(H)},\\
	\operatorname{Var}\bigl(X_{n}(t)\bigl) &= n^{2\mathcal{V}(H)}\left(\mathcal{C}^{\prime}_{n,H}(t)+o(1)\right),
\end{align*}
where the constant $\mathcal{P}_{H}(t)$ is defined in \eqref{P_def}. Hence, it is natural to consider (if it exists) the limit
\begin{align}
	\lim_{n\to\infty}\frac{X_n(t) - \E{X_n(t)}}{n^{\mathcal{V}(H)}}.\label{two_ways_feedback}
\end{align}
Similarly to the computation of $\operatorname{Var}\bigl(X_{n}(t)\bigl)$ done for the one-way feedback model, one can derive higher moments of $X_{n}(t)$, i.e., for any $z\in\N$, $z>2$,
\begin{align*}
	\E{\bigl(X_{n}(t) - \E{X_n(t)}\bigr)^{z}} = n^{z\mathcal{V}(H)}\bigl(\mathcal{C}^{(z)}_{n,H}(t)+o(1)\bigr),
\end{align*}
where the constant $\mathcal{C}^{(z)}_{n,H}(t)$ is defined for $z$ disjoint vectors $\vk v_1,\ldots,\vk v_z\in{\tt K}^{\mathcal{V}(H)}(\scale{n})$ as
\begin{align*}
	\mathcal{C}^{(z)}_{n,H}(t) = \frac{1}{\bigl(\mathcal{V}(H)!\bigr)^z}\sum_{\substack{h_1\in {\tt G}_{\vk v_1}(H) \\\ldots \\ h_{z}\in {\tt G}_{\vk v_z}(H)}}\E{\prod_{k=1}^{z}\left(\mathbb{I}_{n}(h_k,t) - \mathcal{P}_H(t)\right)} = \frac{\E{\prod_{k=1}^{z}\left(\mathbb{I}_{n}(h_k,t) - \mathcal{P}_H(t)\right)}}{\bigl(\mathcal{A}(H)\bigr)^{z}},
\end{align*}
for any choice $h_1\in {\tt G}_{\vk v_1}(H),\ldots,h_z\in {\tt G}_{\vk v_z}(H)$ due to symmetry of the model and \eqref{G_cardinality}.  Hence, if the limit \eqref{two_ways_feedback} exists, its $z$-th moment equals $\lim_{n\to\infty}\mathcal{C}^{(z)}_{n,H}(t)$. Since the random variables $\mathbb{I}_n(h_k,t)$ are dependent for disjoint graphs $h_k$, it is natural to predict that $\lim_{n\to\infty}\mathcal{C}^{(z)}_{n,H}(t)$ should be different from zero for any $z>2$, similarly to $C^{\prime}_n(t)$. However, the simulations suggest that $\lim_{n\to\infty}\mathcal{C}^{(z)}_{n,H}(t) = 0$ (see Table \ref{tab:simulation} for some numerical results  when $n=100$ and $n=200$). This suggests that establishing such limits rigorously would require substantially different techniques, which fall beyond the scope of this work. Therefore, we do not pursue a deeper analysis of this model and leave it as a direction for future research.

\begin{table}
	\centering
	\begin{tabular}{lcccc}
		\toprule
		\textbf{$t$} & \textbf{ $\mathcal{C}^{\prime}_{100,H}(t)$} & \textbf{ $\mathcal{C}^{(3)}_{100,H}(t)$} & \textbf{ $\mathcal{C}^{\prime}_{200,H}(t)$} & \textbf{ $\mathcal{C}^{(3)}_{200,H}(t)$} \\
		\midrule
		1 & 0.1545 &  0.0068 & 0.1565 & 0.0013 \\
		2 & 0.1675 &  0.0109 & 0.1485 & 0.0020 \\ 
		4 & 0.1565 &  0.0145 & 0.1415 &  0.0025 \\
		\bottomrule
		\smallskip
	\end{tabular}
	\caption{Estimation of $\mathcal{C}^{\prime}_{100,H}(t)$, $\mathcal{C}^{(3)}_{100,H}(t)$, $\mathcal{C}^{\prime}_{200,H}(t)$ and $\mathcal{C}^{(3)}_{200,H}(t)$ when $p_0=0.1$, $\beta=0.66$, $\pi_+=0.8$, $\pi_-=0.2$, $H$ is the voter subgraph composed by an edge connecting two vertices having opinion $+$ and $t=1,2,4$. Vertices initially hold opinion $+$ with probability 0.5. Simulations are based on 100 runs.}
	\label{tab:simulation}
\end{table}

\small
\bibliographystyle{abbrvnat}
\bibliography{reference.bib}

\end{document}